\documentclass[11pt]{amsart}
\usepackage{amsmath,amssymb,amsfonts}
\usepackage{enumerate}
\usepackage{amscd, amssymb, latexsym, amsmath, amscd}
\usepackage[all]{xy}
\usepackage{pb-diagram}
\usepackage{picture}
\usepackage{multicol,lipsum}

\newtheorem{theorem}{Theorem}[section]
\newtheorem{thm}[equation]{Theorem}
\newtheorem{prop}[equation]{Proposition}
\newtheorem{cor}[equation]{Corollary}
\newtheorem{con}[equation]{Conjecture}
\newtheorem{lemma}[equation]{Lemma}

\newtheorem{defn}[equation]{Definition}
\numberwithin{equation}{section}

\newcommand{\Z}{\mathbb Z}
\newcommand{\R}{\mathbb R}
\newcommand{\C}{\mathbb C}

\newcommand{\G}{\mathfrak G}

\newcommand{\A}{\mathbb A}

\def\Hom{{\rm Hom}}
\def\Aut{{\rm Aut}}

\def\G{{\rm G}}
\def\SL{{\rm SL}}
\def\disc{{\rm disc}}

\def\PGSp{{\rm PGSp}}

\def\PGSO{{\rm PGSO}}
\def\Spin{{\rm Spin}}
\def\Sp{{\rm Sp}}
\def\Spin{{\rm Spin}}

\def\PU{{\rm PU}}
\def\U{{\rm U}}

\def\GL{{\rm GL}}
\def\PGL{{\rm PGL}}

\def\SO{{\rm SO}}

\def\Sp{{\rm Sp}}
\def\Mp{{\rm Mp}}

\def\O{{\rm O}}
\def\Sym{{\rm Sym}}

\usepackage[all]{xy}

\def\A{{\mathbb A}}

\def\R{{\mathbb R}}

\def\Z{{\mathbb Z}}

\def\C{{\bf C}}

\def\C{{\mathbb C}}

\def\G{{\mathbb G}}

\begin{document}

\title[Theta Correspondence and Automorphic Descent]{Explicit Constructions of Automorphic Forms: \\Theta Correspondence and Automorphic Descent} 
 \author{Wee Teck Gan}
 

\maketitle

This is an expanded set of notes based on two lectures given by the author at the 2022 IHES summer school on the Langlands program, devoted to explicit techniques for the  construction of automorphic representations. In the first lecture, we discussed the classical technique of theta correspondence. The second lecture  covered more recent techniques such as automorphic descent   and the generalized doubling method.  We shall interpret some of these constructions through the lens of the relative Langlands program.
\vskip 10pt

\section{\bf Introduction}

The theta correspondence has always been somewhat of an enigma in the Langlands program. As a topic in representation theory, it has unquestionable pedigree, 
arising from the mathematical foundations of quantum mechanics, where the oscillator representation first made its appearance, and receiving inspiration from the classical invariant theory of Hilbert and Weyl. Weyl's influential book {\em The Classical Groups} \cite{Weyl39} serves as an  important algebraic precursor, with the double centralizer theorem and the Schur-Weyl duality among the highlights. Indeed,  the title of Howe's foundational paper \cite{Howe89}  on the theory of theta correspondence is ``{\em Transcending classical invariant theory}", suggesting that one is going beyond classical invariant theory, by considering a transcendental version of it (which deals with Schwarz functions instead of polynomial functions).

\vskip 5pt

In the theory of automorphic forms, the value of  theta series  as a useful tool for constructing interesting examples of modular forms was recognized very early on, and important arithemetic applications, such as  the Siegel mass formula, had been actively pursued. In the Corvallis proceedings \cite{Corvallis}, which is the standard reference for the Langlands program, there are 4 articles devoted to the theta correspondence, by Howe \cite{Howe79}, Howe-Piatetski-Shapiro \cite{HPS79}, Gelbart \cite{Gelbart79}  and Rallis \cite{Rallis79}. In particular, the paper \cite{HPS79} of Howe-Piatetski-Shapiro used theta correspondence to produce the first counterexamples to the generalized Ramanujan-Petersson conjecture, by constructing nontempered cuspidal representations on the symplectic group $\Sp_4$; this was considered a surprise and breakthrough in the theory of automorphic forms at that time. This was followed shortly by Waldspurger's beautiful rendition of the Shimura correspondence between integral and half-integral weight modular forms \cite{Wald80}, couched in the language of automorphic representations.
\vskip 5pt

Despite the above, the theory of theta correspondence has continued to occupy an uneasy position in relation to the Langlands program. Perhaps it is because the key player in the theory (the Weil representation) is a representation of  a non-linear-algebraic  group (the metaplectic group), which lies outside the traditional realm of the Langlands program.  Or perhaps, its scope  is considered somewhat limited when compared to the ambitious goals of  the Langlands functoriality principle. 
For whatever reasons, the theta correspondence would not have been considered a basic object of study or a fundamental  tool in the classical Langlands program.
\vskip 5pt

A chief purpose of this article is to explain how, with the advent of the {\em relative Langlands program} (as covered in the lectures of Beuzart-Plessis), the theory of theta correspondence has found its natural place within the Langlands philosophy. This requires one to broaden the scope of the relative Langlands program, so that it deals not just with spherical varieties but with quantizations of (certain) Hamiltonian varieties, as discussed in the lectures of Akshay Venkatesh.
\vskip 5pt

Before entering this somewhat philosophical discussion, it will be pertinent to introduce theta correspondence more formally: introducing the main players,  formulating the main questions in the theory and  discussing  some of the progress that have been achieved since the late 1970's.  In particular, we will highlight certain key principles and constructions that have been prominent in the past decades: the tower philosophy, the doubling construction,  the   conservation relation and dichotomy phenomenon, as well as the Rallis inner product formula. We discuss these principles and results first because they are foundational results which belong properly to the theory of theta correspondence (and so exist independently of the Langlands program). We then highlight the role of A-parameters in the theory of theta correspondence: this can be considered an important first point of contact with the Langlands program, with Adams' conjecture playing a decisive role. These considerations lead to 
some spectacular recent applications of theta correspondence to the Langlands program (in the context of non-quasi-split pure inner forms of classical groups).
In addition, we will mention and point towards some outgrowth of classical theta correspondence that are being developed these days. One of these is the extension of theta correspondence beyond the setting of dual pairs in symplectic groups: the {\em exceptional theta correspondence}. Another is the theory of {\em arithmetic theta correspondence}, which will be discussed in the lectures of Chao Li and Wei Zhang.  

\vskip 5pt

We should perhaps mention some topics that will inadvertently be left out due to constraints of space-time and the lack of expertise of the author:
\vskip 5pt

\begin{itemize}
\item the archimedean theta correspondence will go largely unmentioned;

\item there is a geometric theory of theta correspondence in the context of the geometric Langlands program, developed chiefly by S. Lysenko;

\item the geometric theta lifting of Kudla-Millson, which continues to be a very active field of research today;

\item the singular theta lifting of Borcherds.
\end{itemize}
\vskip 5pt

 The last two sections are devoted to various variants of the automorphic descent construction. We first recall the classical automorphic descent of Ginzburg-Rallis-Soudry \cite{GRS2011} and then a variant due to Jiang-Zhang \cite{JiangZhang2020}. Finally, we introduce  the generalized doubling zeta integral of Cai-Fridberg-Ginzburg-Kaplan \cite{CFGK2019}, which represents the standard L-function of $G \times \GL_k$, with $G$ any classical group. This leads to the double descent construction of Ginzburg-Soudry \cite{GinzburgSoudry2021}, which one might consider the ultimate version of automorphic descent. 

\vskip 5pt

\section{\bf The Main Players}
We begin with a traditional introduction to the basic setup of theta correspondence. The standard references for this material are the monograph \cite{MVW87} of Moeglin-Vigneras-Waldspurger and the lecture notes \cite{Kudla93} of Kudla.

 \vskip 5pt

\subsection{\bf Dual pairs} 
 
Consider
\begin{itemize}
\item a field $F$ of characteristic not 2;

\item  a finite dimensional quadratic space $(V,q)$ over $F$, with quadratic form $q : V \rightarrow F$, associated symmetric bilinear form 
\[ \langle  v_1, v_2\rangle_V = q(v_1+ v_2) - q(v_1) -q(v_2),\]
 and isometry group $\O(V)$. For simplicity, we shall assume that ${\rm disc}(V)$ is trivial in $F^{\times}/F^{\times 2}$;

\item  a finite dimensional symplectic vector space $(W, \langle -, - \rangle_W)$  over $F$, with isometry group $\Sp(W)$.
\end{itemize}

Then the tensor product space $\mathbb{W} = V \otimes W$ inherits a natural symplectic form:
\[  \langle -, - \rangle_{\mathbb{W}} = \langle -, - \rangle_V \otimes \langle-,- \rangle_W, \]
and there is a natural map 
\[  i:  \O(V) \times \Sp(W) \longrightarrow \Sp(V \otimes W). \]
The restriction of $i$ to $\O(V)$ or $\Sp(W)$ is injective. The two resulting subgroups are mutual centralizers of each other in $\Sp(V \otimes W)$. Such a pair of reductive subgroups of an ambient group which are mutual centralizers of each other is called a {\em reductive dual pair}. 
\vskip 5pt

Reductive dual pairs in symplectic groups have been completely classified by Howe and discussed in his Corvallis article \cite{Howe79} (see also \cite{MVW87}).
Besides the case of symplectic-orthogonal groups, the other instances are dual pairs  of unitary groups as well as quaternionic unitary groups. More precisely, 
one may consider a pair $V$ and $W$ of Hermitian and skew-Hermitian spaces over appropriate associative algebras with involutions, which results in a dual pair $\U(V) \times \U(W)$.  We will stick with the symplectic-orthogonal example in our exposition; for the case of unitary dual pairs, the reader can consult \cite{Gan-ICM} or \cite{Gan-AWS}.
\vskip 5pt

\vskip 5pt

\subsection{\bf Metaplectic groups}
We now specialize to the case where $F$ is a local field. 
A fact of life  is that in the theory of theta correapondence, we need to work not just with the symplectic group $\Sp(\mathbb{W}) = \Sp(\mathbb{W}_F)$ but its unique nonlinear double cover:
\[  \begin{CD}
1 @>>> \mu_2 @>>> \Mp(\mathbb{W}) @>>> \Sp(\mathbb{W}) @>>> 1. 
\end{CD} \]
Sometimes it will be convenient to push out this extension via $\mu_2 \hookrightarrow S^1$, where $S^1$ is the circle group. We shall also denote such a pushout by $\Mp(\mathbb{W})$ without further comment.  From the structural point of view, a key property of this metaplectic cover is:
\vskip 5pt

\begin{itemize} 
\item  If two elements of $\Sp(\mathbb{W})$ commute  with each other, then any of their preimages commute in $\Mp(\mathbb{W})$.
\end{itemize}
Hence, the preimages in the metaplectic cover of the two members of a dual pair continue to be a dual pair in $\Mp(\mathbb{W})$.
\vskip 5pt

Suppose now $k$ is a global field with adele ring $\A$ and $\mathbb{W}$ is a symplectic vector space over $k$.  Then one has the local metaplectic group $\Mp(\mathbb{W}_v)$ for each place $v$ of $k$. It is known that for almost all $v$,  the metaplectic cover splits uniquely over a hyperspecial maximal compact subgroup $K_v \subset \Sp(\mathbb{W}_v)$. This allows one to define the restricted direct product
\[  \begin{CD} 
1 @>>> \oplus_v  \mu_2 @>>>\prod'_v \Mp(\mathbb{W}_v)   @>>> \Sp(W_{\A})@>>> 1 
\end{CD} \]
Pushing out this extension via the product map $\oplus_v \mu_2 \rightarrow \mu_2$, one obtains a double cover
\[  \begin{CD}
1 @>>> \mu_2 @>>> \Mp(\mathbb{W}_{\A}) @>>> \Sp(\mathbb{W}_{\A}) @>>> 1 \end{CD} 
\]
of the adelic symplectic group. We call $\Mp(\mathbb{W}_{\A})$ the global metaplectic group. A basic property is:
\vskip 5pt

\begin{itemize}
\item This double cover of $\Sp(\mathbb{W}_{\A})$ splits uniquely over the group $\Sp(W_k)$ of rational points.  
\end{itemize}
This property is crucial as it allows us to consider the space  $\mathcal{A}(\Mp(\mathbb{W}))$ of automorphic forms on $\Mp(\mathbb{W}_{\A})$.

\vskip 5pt

\subsection{\bf Weil representations}
The main reason for considering the metaplectic groups is that, over a local field $F$,  $\Mp(\mathbb{W})$ has a finite family of distinguished smooth genuine representations: the Weil representations $\omega_{\mathbb{W},\psi}$. These Weil representations are parametrized by a nontrivial additive character $\psi: F \rightarrow S^1$, with the provision that 
\[  \omega_{\mathbb{W}, \psi} \cong \omega_{\mathbb{W}, \psi'}  \Longleftrightarrow  \text{$\psi'(x) = \psi(a^2 x)$ for some $a \in F^{\times}$. } \]
They arise as  quantizations of the symplectic vector space $\mathbb{W}$ (a theme we shall return to later on), via the representation theory of Heisenberg groups (in particular the Stone-von-Neumann theorem) and their construction goes hand-in-hand with the construction of the metaplectic groups $\Mp(\mathbb{W})$.
In addition, they are unitarizable and can be characterized as the smallest infinite-dimensional representations (in the sense of Gelfand-Kirillov dimension) of $\Mp(\mathbb{W})$.  
\vskip 5pt

\vskip 5pt

\subsection{\bf Schrodinger model} 
One can describe a  concrete model for the Weil representation $\omega_{\mathbb{W}, \psi}$ as follows. If 
\[  \mathbb{W} = \mathbb{X} \oplus \mathbb{Y} \]
 is a Witt decomposition or polarization, so that each $\mathbb{X}$ and $\mathbb{Y}$ is a maximal isotropic subspace, then  $\omega_{\mathbb{W}, \psi}$  can be realized on the space $\mathcal{S}(\mathbb{Y})$ of Schwarz-Bruhat functions on $\mathbb{Y}$. Let $P(\mathbb{X}) = M(\mathbb{X}) \cdot N(\mathbb{X})$ be the Siegel parabolic subgroup stabilizing $\mathbb{X}$ with Levi subgroup $M(\mathbb{X}) \cong \GL(\mathbb{X})$  preserving the Witt decomposition and unipotent radical
\[  N(\mathbb{X})  = \{ n(B): B \in Sym^2 \mathbb{X} \subset \Hom(\mathbb{Y},\mathbb{X}) \}, \]
 where
\[  n(B) = \left( \begin{array}{cc}
1 & B \\
 0 & 1 \end{array} \right)  \quad \text{(relative to $\mathbb{W} = \mathbb{X} \oplus \mathbb{Y}$.)} \]
 The metaplectic cover (pushed out to $S^1$) splits canonically over $N(\mathbb{X})$ and noncanonically over $M(\mathbb{X})$.
 Now one can write down formulas for the action of elements lying over $g \in \GL(\mathbb{X})$ and $n(B) \in N(\mathbb{X})$:
\[ \begin{cases}
(\omega_{\psi}(g) \phi)(y) =   |\det_{\mathbb{X}}(g)|^{\dim \mathbb{W}/2} \cdot \phi(g^{-1} \cdot y)   \\
(\omega_{\psi}(n(B)) \phi)(y)   = \psi( \frac{1}{2} \cdot \langle By, y \rangle) \cdot \phi(y). \end{cases} \]
 To fully describe the action of $\Mp(\mathbb{W})$ (at least projectively), one needs to give the action of an extra Weyl group element 
\[  w = \left( \begin{array}{cc}
0 & 1 \\
 -1 & 0 \end{array} \right)  \]
 which together with $P(\mathbb{X})$ generates $\Sp(\mathbb{W})$. The action of this $w$ is given by Fourier transform  (up to scalars):
 \[  (\omega_{\psi}(w) \phi)(y) = \int_{\mathbb{Y}} \phi(z) \cdot  \psi \left( \langle w \cdot z, y  \rangle \right) \, dz  \]
 So we see  the ubiquitous tool of Fourier transform  contained within the theory of the Weil representation.
\vskip 5pt

\vskip 5pt
\subsection{\bf Theta functions}
Over a global field $k$, for each nontrivial $\psi: k \backslash \A \rightarrow S^1$, one  has the abstract Weil representation
 \[  \omega_{\mathbb{W}_{\A}, \psi}  = \bigotimes'_v \omega_{\mathbb{W}_v, \psi_v} \]
of the adelic metaplectic group $\Mp(\mathbb{W}_{\A})$, which can be realized on the space 
\[  \mathcal{S}(\mathbb{Y}_{\A}) = \otimes'_v \mathcal{S}(\mathbb{Y}_v). \]
The key fact about this abstract adelic representation is that it admits a natural equivariant map
\[  \theta:  \mathcal{S}(\mathbb{Y}_{\A}) \longrightarrow \mathcal{A}(\Mp(\mathbb{W})) \]
 given by ``averaging over rational points":
 \[  \theta(\phi)(g)   =   \sum_{y \in \mathbb{Y}_k} (\omega_{\psi}(g) \phi)(y). \]
The fact that the function $\theta(\phi)$ is left-invariant under $P(\mathbb{X}_k)$ is easy to check using the above formulas in the Schrodinger model. The fact that it is
left-invariant under the Weyl group element $w$ is a consequence of the Poisson summation formula.
The  function $\theta(\phi)$ is the automorphic incarnation of what-are-classically-called theta functions. Hence the map $\theta$ is also called the ``formation of theta series".
\vskip 5pt

\subsection{\bf Splitting of metaplectic covers}
Suppose one has a dual pair 
\[  i: \O(V) \times \Sp(W) \longrightarrow \Sp(\mathbb{W}). \]
To pullback a Weil representation to $\O(V) \times \Sp(W)$, one needs to address the question of whether
 the morphism $i$ can be lifted to $\Mp(\mathbb{W})$ (both locally and globally).  This is a technical issue  and we summarize the results here:
\vskip 5pt

\begin{itemize}
\item If $\dim V$ is even, then the map $i$ can be lifted to a map to $\Mp(\mathbb{W})$ (both locally and globally);

\item if $\dim V$ is odd, then $i|_{\O(V)}$ can be lifted to $\Mp(\mathbb{W})$ but not $i|_{\Sp(W)}$.  
\end{itemize}
In the following, we shall assume implicitly (unless otherwise stated) that $\dim V$ is even, 
so that we may work with the more familiar $\Sp(W)$ instead of $\Mp(W)$.
\vskip 5pt

Note that while the lifting of $i$ is unique on $\Sp(W)$ (if it exists), it is not unique on $\O(V)$, as one can twist any such splitting by a quadratic character of $\O(V)$. A careful construction and enumeration of such splittings has been provided by Kudla \cite{Kudla94}. Without going into details, one would say that the data of $(W, \psi)$ (with $\psi$ a nontrivial additive character of $F$ in the local setting and a character $\psi: k \backslash \A$ in the global setting) defines a lifting
\[  i_{W,\psi}: \O(V) \longrightarrow  \Mp(\mathbb{W}). \]
One can describe this splitting concretely using the Schrodinger model.
If $W = X \oplus Y$ is a polarization, so that $\mathbb{Y} = V \otimes Y$ is a maximal isotropic subspace of $\mathbb{W}$, then one has the Schrodinger model of $\omega_{\mathbb{W}, \psi}$ realized on $\mathcal{S}(\mathbb{Y})$. In this model, the action of $\iota_{W,\psi}(\O(V))$ is geometric:
\[  (h \cdot \phi) (y) = \phi( h^{-1} \cdot y), \quad \text{  for $h \in \O(V)$.} \]
\vskip 5pt

In any case, by using a lifting (whenever it exists)
\[ i_{V, W,\psi} : \O(V) \times \Sp(W) \longrightarrow \Mp(\mathbb{W}), \]
one may pull back the Weil representation $\omega_{\psi}$ to obtain a representation
\[  \Omega_{V,W, \psi} = \omega_{\psi} \circ i_{V,W,\psi} \quad \text{  of $\O(V) \times \Sp(W)$} \]
in the local setting.  
\vskip 5pt

Likewise, in the global setting, one obtains by restriction of functions an equivariant map
\[  \begin{CD}
\Omega_{V,W, \psi} = \mathcal{S}(\mathbb{Y}_{\A})  @>\theta>>   \mathcal{A}(\Mp(\mathbb{W})) @>i_{V,W,\psi}^*>>  \mathcal{C}([\O(V) \times \Sp(W)]) \end{CD} \]
where the target is the space of smooth functions on 
\[  [\O(V)] \times [\Sp(W)] = \O(V_k) \backslash \O(V_{\A})  \times \Sp(W_k) \backslash \Sp(W_{\A}). \]

\vskip 10pt
\section{\bf The Problems}
After introducing the main players, we can now formulate the main problems in the theory of theta correspondence.
\vskip 5pt

\subsection{\bf Local restriction problems}
 
In the local setting,  one studies   the spectral decomposition of the representation $\Omega := \Omega_{V, W, \psi}$ as a representation of $\O(V) \times \Sp(W)$. We can ask this question in two different settings:
\vskip 5pt

\begin{itemize}
  
\item (Smooth setting)
Writing ${\rm Irr}(\O(V))$ for the set of equivalence classes of irreducible smooth representations of $\O(V)$ (and likewise  ${\rm Irr}(\Sp(W))$), one is interested in determining all $(\pi, \sigma) \in {\rm Irr}(\O(V)) \times {\rm Irr}(\Sp(W))$ such that
\[  \Hom_{\O(V) \times \Sp(W)} (\Omega_{V,W,\psi}, \pi \otimes \sigma) \ne 0. \]
One can reformulate the above definition in a slightly different way. 
  For $\pi \in {\rm Irr}(\O(V))$, one considers the maximal $\pi$-isotypic quotient of $\Omega$:
\[  \Omega /  \bigcap_{f \in \Hom_{\O(V)}(\Omega, \pi)}  {\rm Ker} (f),  \]
which is a $\O(V) \times \Sp(W)$-quotient of $\Omega$ expressible in the form
\[ \pi \otimes \Theta(\pi), \]
for some smooth representation $\Theta(\pi)$ of $\Sp(W)$ (possibly zero, and possibly infinite length a priori).
 We call $\Theta(\pi)$ the {\em big theta lift} of $\pi$. 
 A more direct way to define $\Theta(\pi)$ is:
\[  \Theta(\pi) = ( \Omega \otimes \pi^{\vee})_{\O(V)}, \]
 the maximal $\O(V)$-invariant quotient of $\Omega \otimes \pi^{\vee}$. It follows from definition that there is a natural $\O(V) \times \Sp(W)$-equivariant map
 \[  \Omega \twoheadrightarrow \pi \otimes \Theta(\pi), \]
 which satisfies the ``universal property":  for any smooth representation $\sigma$ of $\Sp(W)$, 
 \[  \Hom_{\O(V) \times \Sp(W)}(\Omega, \pi \otimes \sigma) \cong \Hom_{\Sp(W)}(\Theta(\pi), \sigma)  \quad \text{(functorially)}. \] 
 The goal of local theta correspondence is to determine the representation $\Theta(\pi)$ or rather its irreducible quotients. By symmetry, one also has the representation $\Theta(\sigma)$ of $\O(V)$ for $\sigma \in {\rm Irr}(\Sp(W))$.
   \vskip 10pt

\item ($L^2$-setting)
  We may also consider this restriction problem in the $L^2$-setting. Recall that the Weil representation $\Omega$ is unitairizable and hence we may consider its unitary completion. In the Schrodinger model, this would be a unitary representation realized on the Hilbert space $L^2(\mathbb{Y})$.  We denote this unitary Weil representation by $\hat{\Omega}$. Then one has a direct integral decomposition
  \[  \hat{\Omega} \cong \int_{\widehat{\O(V)}} \hat{\pi} \otimes \hat{\theta}(\hat{\pi})  \, d\mu(\hat{\pi}) \]
for some unitary representation $\hat{\theta}(\hat{\pi})$ of $\Sp(W)$ and for some measure $d\mu$ on the unitary dual $\widehat{\O(V)}$ of $\O(V)$. 
Here the main issue is to determine the measure $d\mu$ and the representation $\hat{\theta}(\hat{\pi})$ for $d\mu$-almost all $\hat{\pi} \in \widehat{\O(V)}$. A more precise problem would be to explicate the unitary isomorphism between the two sides. One way of doing so is to give a spectral decomposition of the $L^2$-inner product of $\hat{\Omega}_{V,W,\psi}$, i.e.
\[  \langle \phi_1, \phi_2 \rangle  =  \int_{\widehat{\O(V)}}  J_{\hat{\pi}}(\phi_1, \phi_2) \, d\mu(\hat{\pi}) \]
where $J_{\hat{\pi}}$ is the $\hat{\pi} \otimes \hat{\theta}(\hat{\pi})$-component of the inner product $\langle-, -\rangle$: it is a degenerate equivariant pairing on $\hat{\pi} \otimes \hat{\theta}(\hat{\pi})$ which one would like to explicate.
 
\end{itemize}
\vskip 5pt

We note that the smooth and $L^2$-version of the local theta correspondence are related. More precisely, for $d\mu$-almost all $\hat{\pi}$, the spectral decomposition of $\hat{\Omega}$ gives rise to an equivairant projection
\[  \Omega= \mathcal{S}(\mathbb{Y})  \longrightarrow \pi \otimes \hat{\theta}(\hat{\pi})^{\infty}  \]
where $\pi = \hat{\pi}^{\infty}$ and $\hat{\theta}(\hat{\pi})^{\infty}$ denote the underlying smooth representations. Hence, by the ``universal property" alluded to in the smooth setting, $\hat{\theta}(\hat{\pi})^{\infty}$ is a semisimple quotient of the smooth representation $\Theta(\pi)$.
\vskip 10pt

\subsection{\bf Global theta liftings}  \label{SS:global theta}
Let us now consider the global setting where we work over a number field $k$ with adele ring $\A$. 
We have noted that the adelic Weil representation  $\Omega_{V,W,\psi}$ of the dual pair $\O(V) \times \Sp(W)$ admits an equivariant map
\[  \begin{CD}
\Omega_{V,W, \psi} = \mathcal{S}(\mathbb{Y}_{\A})  @>\theta>>   \mathcal{A}(\Mp(\mathbb{W})) @>i_{V,W,\psi}^*>>  \mathcal{C}([\O(V) \times \Sp(W)]) \end{CD} \]
given by the ``formation of theta series".  For each $\phi \in \Omega_{V,W, \psi}$, $\theta(\phi)$ is thus a function on $[\O(V) \times \Sp(W)]$ are thus can be used as a kernel function to transfer functions on $[\O(V)]$ to those on $[\Sp(W)]$. 
\vskip 5pt

More precisely, suppose that $\sigma \subset \mathcal{A}_{cusp}(\Sp(W))$ is a cuspidal representation of $\Sp(W)(\A)$. Then for $\phi \in \Omega_{V,W,\psi}$ and $f \in \sigma$, we set:
\[  \theta(\phi,f)(h) := \int_{[\Sp(W)]} \theta(\phi) (gh) \cdot \overline{f(g)} \, dg, \]
for $g \in \Sp(W_{\A})$ and $dg$ denoting the Tamagawa measure.  Then we set
\[  \Theta(\sigma) := \langle \theta(\phi,f): \phi \in \Omega_{V,W,\psi}, \, f \in \sigma \rangle \subset \mathcal{A}(\O(V)). \]
This is a $\O(V_{\A})$-submodule of the space of automorphic forms on $\O(V)$ and we call it the global theta lift of $\sigma$.
\vskip 5pt

The basic questions in global theta correspondence parallel those in the local setting:
\vskip 5pt

\begin{itemize}
\item Decide if $\Theta(\sigma)$ is nonzero.
\item Is  $\Theta(\sigma) \subset \mathcal{A}_{cusp}(\O(V))$? Is $\Theta(\sigma) \subset \mathcal{A}_2(\O(V))$?
\item What automorphic representations do one get from this construction?
\item What is the relation between the global theta lift and the  local ones?
\end{itemize}

 \vskip 10pt

\section{\bf Some Local Results}
Having defined the problems studied in the theta correspondence, we can begin to formulate some results in the local setting, so that we are working over a local field $F$. As mentioned in the introductory subsection, we will focus first on foundational results within the theory proper which are independent of the Langlands program. 
\vskip 5pt

\subsection{\bf Howe duality theorem}
We  first note the following basic result of Howe \cite{Howe89} and Kudla \cite{Kudla86}:
\begin{prop}[Finiteness]
For any $\pi \in {\rm Irr}(\O(V))$, $\Theta(\pi)$ is of finite length. In particular, if $\Theta(\pi)$ is nonzero, it has (finitely many) irreducible quotients and we may consider its maximal semisimple quotient (its cosocle) $\theta(\pi)$. Moreover, for any $\pi \in{\rm Irr}(\O(V))$ and $\sigma \in {\rm Irr}(\Sp(W))$,
\[  \dim \Hom_{\O(V) \times \Sp(W)}(\Omega, \pi \otimes \sigma) < \infty. \]
 \vskip 5pt
 
 \end{prop}
We call $\theta(\pi)$ the {\em small theta lift} of $\pi$. The local theta lifts of $\pi$ are precisely the irreducible summands of $\theta(\pi)$.  
 \vskip 5pt

We can now formulate a fundamental theorem which refines the above finiteness result. This so-called Howe duality theorem was shown by Howe \cite{Howe89} in the archimedean case, by Waldspurger \cite{Wald90} in the p-adic case when $p \ne 2$ and in general in \cite{Minguez2008, GanTakeda2016, GanSun2017}.
\vskip 5pt

\begin{thm}[Howe duality theorem]
(i) If $\Theta(\pi)$ is nonzero, then it has a unique irreducible quotient. In other words,  $\theta(\pi)$ is irreducible or $0$.
\vskip 5pt

(ii) If $\theta(\pi) \cong \theta(\pi') \ne 0$, then $\pi \cong \pi'$.
 
\vskip 5pt

In other words, one has:
\[  \dim \Hom_{\Sp(W)}( \theta(\pi) , \theta(\pi')) \leq \dim \Hom_{\O(V)}(\pi, \pi') \leq 1, \]
so that one has a map
\[  \theta_{V,W, \psi} : {\rm Irr}(\O(V)) \longrightarrow  {\rm Irr}(\Sp(W)) \cup \{ 0 \} \]
which is injective when restricted to the subset of ${\rm Irr}(\O(V))$ which is not sent to $0$. 
 
\end{thm}
\vskip 5pt
  Of course, by symmetry, we have the analogous theorem with the roles of $\O(V)$ and $\Sp(W)$ exchanged.  
  
 \vskip 10pt

\subsection{\bf Further Questions}  \label{SS:further questions}
After the Howe duality theorem above, we can ask the following more refined questions:
\vskip 5pt

\begin{itemize}
\item[(a)] (Nonvanishing)  For a given $\pi \in {\rm Irr}(\O(V))$, decide if $\theta(\pi)$ is nonzero.
\vskip 5pt

\item[(b)] (Identity) Describe the map $\theta$ explicitly. In other words,  If $\theta(\pi)$ is nonzero, can one describe it in terms of $\pi$ in terms of, for example, the local Langlands correspondence (or some other way).    
\end{itemize}
   \vskip 5pt

\subsection{\bf Rallis' tower}
For the nonvanishing question,  Rallis observed that it is fruitful to consider theta correspondence in a family. 
Let  
\[  W_r = \mathbb{H}^r \]
be the $2r$-dimensional symplectic space,  where $\mathbb{H}$ is the hyperbolic plane (i.e. the 2-dimensional symplectic space).
The collection $ \mathcal{W} = \{ W_r \, | \, r \ge 0 \}$ is called a Witt tower of symplectic spaces. 
One can likewise consider a Witt tower of quadratic spaces. Fixing an anisotropic quadratic space $V_0$, we se
\[  V_r = V_0 \oplus \mathbb{H}^r \]
where $\mathbb{H}$ is the hyperbolic plane (i.e. the split 2-dimensional quadratic space).
Observe that:
\begin{itemize}
\item $\dim V_r \mod 2$ is independent of $r$;
\item  $\disc(V_r) \in F^{\times}/F^{\times 2}$  is independent of $r$.
\end{itemize}
Hence, in the nonarchimedean case, there are precisely two Witt towers of quadratic spaces with a fixed dimension modulo 2 and fixed discriminant in $F^{\times}/ F^{\times 2}$, which are distinguished by their Hasse-Witt invariant $\epsilon(V) = \pm 1$. We shall call these the plus- and minus-tower respectively, denoting them by $\mathcal{V}^+ = \{V_r^+ \}$ and $\mathcal{V}^- = \{ V_r^- \}$ respectively. Any symplectic or quadratic space is a member of a unique Witt tower.

 \vskip 5pt
If we fix a sympelctic space $W$,  we may  consider a family of theta correspondences associated to the tower of reductive dual pairs:
\vskip 5pt

\begin{itemize}
\item  For a fixed symplectic space $W$ and $\sigma \in {\rm Irr}(\Sp(W))$, we may consider its theta lift $\theta^{\epsilon}_r(\sigma)$ on $\O(V^{\epsilon}_r)$, relative to the tower of Weil representations $\Omega_{W, \mathcal{V}^{\epsilon},\psi}$.
\vskip 5pt

\item For a fixed quadratic space $V$ and $\pi \in {\rm Irr}(\O(V))$, one may consider its theta lift $\theta^+_r(\pi)$ on $\Sp(W_r)$. 
\end{itemize}
\vskip 5pt

\noindent It would appear that there is an asymmetry between the symplectic and orthogonal case, since there is only one Witt tower of symplectic spaces and hence only one Witt tower of theta lifts. However, one should regard a Rallis tower as a tower of representations of the fixed group $\O(V)$, rather than  a tower of spaces (of opposite parity).
Treating the Rallis tower as the plus-tower of Weil representations $\Omega_{V, \mathcal{W}, \psi}$, we may consider also the complementary tower of twisted Weil representations $\Omega_{V, \mathcal{W}, \psi} \otimes {\det}_V$ and call this the minus-tower of Weil representation. Thus, we may set
\[  \theta^-_r(\pi) := \theta^+_r(\pi \otimes \det). \]
For a more uniform (or symmetrical) formulation, the reader can consult \cite{SunZhu2015}. 
\vskip 5pt

Now Kudla \cite{Kudla86, Kudla93} showed:
\vskip 5pt

\begin{prop} \label{P:kudla86}
(i) (First occurrence) For fixed $\epsilon = \pm $ and $\pi \in  {\rm Irr}(\O(V))$, there is a smallest $r^{\epsilon}_0 = r^{\epsilon}_0(\pi)$ such that $\theta^{\epsilon}_{r^{\epsilon}_0}(\pi)\ne 0$. Moreover, $r^{\epsilon}_0 \leq \dim V$. This $r^{\epsilon}_0(\pi)$ is called the first occurrence index  of $\pi$ in the Witt tower $\mathcal{W}^{\epsilon}$.
\vskip 5pt

(ii) (Tower property) For any $r > r^{\epsilon}_0$,   $\theta^{\epsilon}_r (\pi)\ne 0$.  
\vskip 5pt

(iii) (Supercuspidal case) Suppose that $\pi$ is supercuspidal. Then $\Theta_{r^{\epsilon}_0}^{\epsilon}(\pi)$ is irreducible supercuspidal. For $r > r^{\epsilon}_0$, 
$\Theta_r^{\epsilon}(\pi)$ is irreducible but not supercuspidal.
\vskip 5pt

(iv) (Cuspidal support). In general,  the map $\pi \mapsto \Theta_r(\pi)$ commutes with the consideration of cuspidal support. In particular, all irreducible subquotients of $\Theta_r(\pi)$ have the same cuspidal support.

\end{prop}
One has the analogous result starting with a $\sigma \in {\rm Irr}(\Sp(W))$.  Some remarks:

 \vskip 5pt
 
 \begin{itemize}
  
 \item The fact that $r^{\epsilon}_0(\pi)  \leq  \dim V$ means that when $W$ is sufficiently large (more precisely when $W$ has an isotropic  subspace of dimension $\geq \dim V$), the map $\theta_{V,W, \psi}$ is nonzero on the whole of ${\rm Irr}(\O(V))$. When $r \geq \dim V$, we say that the theta lifting is in the {\em stable range}.

\item The nonvanishing problem (a) highlighted above is reduced to the question of determining the first occurrence indices for all the Witt towers.
\end{itemize}
 \vskip 5pt
 
 \subsection{\bf Conservation and Dichotomy}
In fact, our job is halved because the two first occurrence indices  $r_0^{\pm}(\pi)$ are not independent of each other. Rather, we have the following theorem of Sun-Zhu \cite{SunZhu2015} confirming a conjecure of Kudla-Rallis:
\vskip 5pt

\begin{thm}[Conservation relation]
(i) For $\pi \in {\rm Irr}(\O(V))$, let $r^+_0(\pi)$ and $r^-_0(\pi)$ be the two first occurrence indices in the two symplectic towers defined above. Then
\[  r^+_0(\pi) + r^-_0(\pi) = \dim V. \]
\vskip 5pt

(ii) For $\sigma \in {\rm Irr}(\Sp(W))$, let $r^+(\sigma)$ and $r^-(\sigma)$ be the two first occurrence indices in the two quadratic towers (with fixed dimension modulo $2$ and fixed discriminant). Then
\[   r^+_0(\sigma) + r^-_0(\sigma) = \dim W. \] 
\end{thm}
The  conservation relation above implies the following dichotomy theorem:

\begin{cor}[Dichotomy]
Suppose that $W^+$ and $W^-$ belong to the plus- and minus-Witt towers respectively, and $\dim W^+ + \dim W^- = 2 \dim V$. 
Then for any $\pi \in {\rm Irr}(\O(V))$,  exactly one of the two theta lifts $\Theta_{V,W^+,\psi}(\pi)$ and $\Theta_{V, W^-, \psi}(\pi)$ is nonzero.
\end{cor}
Such a dichotomy result was first shown in the unitary case in \cite{HKS96}.
\vskip 5pt

\subsection{\bf Big Theta versus small theta}   \label{SS:big-v-small}
In Proposition \ref{P:kudla86}(iii) above, one sees that for supercuspidal $\pi$, one has $\Theta(\pi) = \theta(\pi)$.  This question of whether the big theta lift of $\pi$ is  equal to the small theta lift of $\pi$ is quite important, as it intervenes in applications of theta correspondence to branching problems (a theme we shall return to later). What do we know beyond the supercuspidsl case? 
\vskip 5pt

In a series of papers \cite{Muic2004, Muic2006, Muic2008a, Muic2008b}, Mui\'c has studied the structure of $\Theta(\pi)$ for discrete series representations $\pi$ and more generally for tempered $\pi$.   We summarize some of his results  and the results of some others \cite{GanSavin2012, GanIchino2014, LokeMa} here:
\vskip 5pt

\begin{thm}
(i) Suppose that $\pi \in {\rm Irr}(\O(V))$ is a discrete series representation and $\dim W \leq \dim V$. Then $\Theta(\pi) = \theta(\pi) \in {\rm Irr}(\Sp(W))$.
The same holds for tempered $\pi$ if $\dim W = \dim V$.
\vskip 5pt

(ii) In the archimedean case, if $\dim W \geq 2 \dim V$ (i.e. in the stable range) and $\pi$ is unitary, then $\Theta(\pi) = \theta(\pi)$.
\end{thm}
Here, part (ii) of the theorem is a result of Loke-Ma \cite{LokeMa}. The result in (i) plays an important role in applications of theta correspondence to the Gross-Prasad conjecture \cite{GanIchino2016}.

\vskip 10pt

 \section{\bf A Basic Tool: the Doubling Seesaw}
 There are two basic ingredients involved in the proof of many of the results mentioned in the previous section:
 \vskip 5pt
 \begin{itemize}
 \item a computation of the Jacquet modules of the Weil representation $\Omega$ relative to a parabolic subgroup of one member of a dual pair; in the local setting, this computation was carried out by Kudla \cite{Kudla86} and in the global setting, the analogous constant term computation was carried out by Rallis.
 \vskip 5pt
 \item the theory of the doubling seesaw.
  \end{itemize}
 
\noindent  We  will discuss this doubling seesaw in this section.
\vskip 5pt

\subsection{\bf Seesaw pairs and identities}  \label{SS:seesaw}
  Let us first introduce a general concept introduced by Kudla \cite{Kudla84}, which he called ``a seesaw pair". This consists of two dual pairs $(G, H)$ and $(G', H')$ in a fixed $\Sp(\mathbb{W})$ such that $G \subset G'$ (and hence $H' \subset H$). This situation can be represented pictorially as:

\[
    \xymatrix{
    G' \ar@{-}[d]\ar@{-}[dr]&  H\ar@{-}[d]\\
    G \ar@{-}[ur]& H'
    }
\]
so that the groups on each diagonal form a dual pair and each group in the bottom row is contained in the group directly above it. This pictorial presentation explains  the terminology ``seesaw pair" . 
 \vskip 5pt
 
 When one has such a seesaw pair, one has an associated  ``seesaw identity". More precisely, given $\pi \in {\rm Irr}(G)$ and $\sigma \in {\rm Irr}(H')$ (i.e. irreducible representations of the groups in the bottom row), one may consider
 \[  \Hom_{G \times H'}(\Omega,  \pi \otimes \sigma) \]
 and compute it in two different ways:
 \[  \Hom_{H'}(\Theta(\pi), \sigma) \cong \Hom_{G \times H'}(\Omega,  \pi \otimes \sigma) \cong \Hom_G(\Theta(\sigma), \pi). \]
 This is a triviality, but a surprisingly useful one, as it allows one to transfer a branching problem from the $G$-side to one on the $H$-side,  often offering a new perspective on the problem.
 \vskip 5pt

 \subsection{\bf The doubling seesaw}  \label{SS:doubling seesaw}
 The seesaw which is most relevant to our discussion is  the so-called doubling see-saw:
\[
    \xymatrix{
    \Sp(W \oplus W^-)\ar@{-}[d]\ar@{-}[dr]& \O(V)\times \O(V)\ar@{-}[d]\\
    \Sp(W)\times \Sp(W^-)\ar@{-}[ur]& \O(V)^{\Delta}
    }
\]
where $W^-$ denotes the symplectic space  obtained from $W$ by multiplying the form by $-1$ (so that $\Sp(W^-) = \Sp(W) \subset \GL(W)$). Hence, we have doubled the space $W$ to obtain $W \oplus W^-$. Note that the diagonally embedded $W^{\Delta}$ is a maximal isotropic subspace of the doubled space $W \oplus W^-$ and one has a Witt decomposition
\[  W + W^- = W^{\Delta} \oplus W^{\nabla} \]
with $W^{\nabla} = \{ (w, -w): w \in W \}$. 
The stabilizer of $W^{\Delta}$ in $\Sp(W+W^-)$ is thus a Siegel maximal parabolic subgroup $P(W^{\Delta})$. We stress that this seesaw already appeared in Howe's Corvallis article. 
\vskip 5pt

Let us apply the seesaw identity for this doubling seesaw to the representations:
\[  \sigma' \otimes \sigma^{\vee} \in {\rm Irr}(\Sp(W) \times \Sp(W')) \quad \text{and} \quad  1 \in {\rm Irr}(\O(V)^{\Delta}). \]
 The seesaw identity gives:
\[  \Hom_{\O(V)}(\Theta_{V,W, \psi}(\sigma') \otimes \Theta_{V,W^-, \psi}(\sigma^{\vee}), \C)   \cong \Hom_{\Sp(W) \times \Sp(W^-)} ( \Theta(1),  \sigma' \otimes \sigma^{\vee}). \]

\vskip 5pt

\subsection{\bf Sketch proof of Howe duality}  \label{SS:proof of Howe}
 The relevance of the above  identity to the Howe duality theorem comes from the observation that
 \[  \Hom_{\O(V)}(\theta_{V,W, \psi}(\sigma'), \theta_{V,W,\psi}(\sigma)) \subset  \Hom_{\O(V)}(\Theta_{V,W, \psi}(\sigma') \otimes \Theta_{V,W^-, \psi}(\sigma^{\vee}), \C). \]
 This requires relating
 \[  \Theta_{V,W^-, \psi}(\sigma^{\vee}) \quad \text{and} \quad \Theta_{V,W, \psi}(\sigma), \]
 a detail we won't go into here (but which involves the so-called MVW involution).  One would thus like to show  that
 \[  \dim \Hom_{\Sp(W) \times \Sp(W^-)} ( \Theta(1),  \sigma' \otimes \sigma^{\vee}) \leq \dim \Hom_{\Sp(W)} ( \sigma' ,  \sigma). \]
To proceed further, one needs to understand the big theta lift $\Theta(1)$ of the trivial representation.  
\vskip 5pt

\subsection{\bf Understanding $\Theta(1)$} 
  Rallis \cite{Rallis84a} showed that the big theta lift $\Theta(1)$ of the trivial representation embeds naturally into a member of the family of Siegel degenerate principal series 
of $\Sp(W + W^-)$  associated to the Siegel parabolic subgroup $P(W^{\Delta})$  
\[    I_{P(W^{\Delta})} (s) = {\rm Ind}_{P(W^{\Delta})}^{\Sp(W+W^-)} |\det|^{s}. \]
To define this embedding, observe that we may realize the Weil representation of $\Sp(W + W^-) \times \O(V)$ in its Schrodiger model $\mathcal{S}(V \otimes W^{\nabla})$. 
Then Rallis considered the natural map
\[   \mathcal{S}(V \otimes W^{\nabla}) \longrightarrow I_{P(W^{\Delta})}(s_{V,W}), \quad \text{with $s_{V,W} = \frac{\dim V - \dim W -1}{2}$,}  \]
defined by 
\[   f_{\phi} (g ) =  (\omega(g) \cdot \phi)(0) \quad\text{  for $\phi \in  \mathcal{S}(V \otimes W^{\nabla})$.} \]
This map is clearly $\O(V)$-invariant and Rallis \cite[Thm. II.1.1]{Rallis84a} showed that it gives an injection
\[  \Theta(1)  \hookrightarrow I_{P(W^{\Delta})}(s_{V,W}). \]
If we assume without loss of generality that $\dim W \geq \dim V- 1$ (so that $s_{V,W} \leq  0$), then one can show that $\Theta(1)$ is irreducible and is a quotient of $I_{P(W^{\Delta})}(-s_{V,W})$.  Thus, one has
\[  
\dim  \Hom_{\O(V)}(\theta(\sigma'), \theta(\sigma))  \leq  \dim \Hom_{\Sp(W) \times \Sp(W^-)} (I_{P(W^{\Delta})} (-s_{V,W}) ,  \sigma' \otimes \sigma^{\vee}). \] 
Hence, one needs to understand $I_{P(W^{\Delta})} (s_{V,W})$ as a $\Sp(W) \times \Sp(W^-)$-module. This is a question of Mackey theory. 

\vskip 10pt

\subsection{\bf Mackey theory}
It turns out that there is a unique open dense $\Sp(W) \times \Sp(W^-)$-orbit on the partial flag variety $P(W^{\Delta}) \backslash \Sp(W + W^-)$, namely the orbit of the maximal isotropic subspace $W^{\Delta}$. Moreover, the stabilizer of $W^{\Delta}$  is the diagonally embedded
\[  \Sp(W)^{\Delta}  \hookrightarrow \Sp(W) \times \Sp(W^-). \]
Thus, we have
\[  C^{\infty}_c(\Sp(W)) \subset I_{P(W^{\Delta})} (-s_{V,W})  \]
as a $\Sp(W) \times \Sp(W-)$-submodule and the quotient is small (since it corresponds to the complement of the open orbit).  
So, up to a first approximation, the natural restriction map
\[   \Hom_{\Sp(W) \times \Sp(W^-)} (I_{P(W^{\Delta})} (-s_{V,W}) ,  \sigma' \otimes \sigma^{\vee}) \longrightarrow \Hom_{\Sp(W) \times \Sp(W^-)} (C^{\infty}_c(\Sp(W) ,  \sigma' \otimes \sigma^{\vee}) \]
is an isomorphism for most $\sigma' \otimes \sigma$. Of course, the latter space is at most 1-dimensional, with equality if and only if $\sigma' \cong \sigma$. 
\vskip 5pt
Hence, we see that the doubling seesaw argument takes care of the Howe duality theorem for almost all representations. The remaining representations (called representations on the boundary) have to be dealt with by some other means. 
\vskip 5pt

\subsection{\bf Dichotomy in the equal rank case}  \label{SS:dichotomy}
  The doubling seesaw also features heavily  in the proof of the conservation relation and the dichotomy theorem. We illustrate this by giving a sketch of dichotomy in the equal rank case, i.e. for the dual pair $\O(V) \times \Mp(W)$ with
 \[  \dim W = 2n \quad \text{and} \quad \dim V = 2n+1. \]
In the seesaw argument above, one has $s_{V,W} = 0$ and 
\[  \Theta_{V, W+W^-}(1) \subset  I_{P(W^{\Delta})}(0). \]
Indeed, one has
\[  I_{P(W^{\Delta})}(0)  = \bigoplus_{\epsilon= \pm} \Theta_{V^{\epsilon}, W+W^-}(1). \] 
The doubling seesaw identity (applied with $\sigma' = \sigma$)  thus gives:
\[  \sum_{\epsilon = \pm}  \dim \Hom_{\O(V)} (\theta_{W,V^{\epsilon}, \psi}(\sigma), \theta_{W,V^{\epsilon}, \psi}(\sigma)) = \dim \Hom_{\Mp(W) \times \Mp(W^-)} (I_{P(W^{\Delta})}(0), \sigma \otimes \sigma^{\vee}). \]
In this case, one can show that tempered representations do not live on the boundary, so that  
\[  \dim \Hom_{\Mp(W) \times \Mp(W^-)} (I_{P(W^{\Delta})}(0), \sigma \otimes \sigma^{\vee}) = \dim \Hom_{\Mp(W) \times \Mp(W^-)} (C^{\infty}_c(\Mp(W)), \sigma \otimes \sigma^{\vee}) = 1\]
for tempered $\sigma$.
This shows that exactly one of 
\[   \theta_{W, V^+,\psi}(\sigma) \quad \text{or} \quad \theta_{W, V^-,\psi}(\sigma)   \]
is nonzero, thus showing dichotomy for tempered representations. 
\vskip 5pt

Of course, this leads to the natural question of which one of the above two theta lifts is nonzero. This turns out to have  an  interesting answer \cite{HKS96, GanSavin2012}:  
\vskip 5pt

\begin{thm}[Epsilon dichotomy]
 Let $\sigma$ be an irreducible representation of $\Mp(W)$. Then $\theta_{W,V^{\epsilon}, \psi}(\sigma) \ne 0$ if and only if 
\[  \epsilon = z_{\psi}(\sigma)     \cdot 
\epsilon(1/2, \sigma, \psi), \]
where $z_{\psi}(\sigma) $ is a sign encoding the central character of $\sigma$  and the local epsilon factor is the standard epsilon factor defined by the doubling zeta integral.
\end{thm}
\vskip 5pt

\section{\bf The Doubling Zeta Integral}
The above theorem compels us to introduce the important doubling zeta integrals. This is a family of Rankin-Selberg integrals discovered by Piatetski-Shapiro and Rallis \cite{GPSR87} and gives an integral representation of the standard L-function of $G \times \GL_1$, where $G$ is any classical group. We will focus on the local theory in this section.
\vskip 5pt

\subsection{\bf The local doubling zeta integral}  \label{SS:local doubling}
  Recall that
\[  \Sp(W^{\Delta}) \backslash \Sp(W) \times \Sp(W^-) \cong \Sp(W)  \hookrightarrow P(W^{\Delta})  \backslash \Sp(W+W^-) \]
as an open dense subset. One thus has a restriction map
\[  I_{P(W^{\Delta})}(s, \chi) := {\rm Ind}_{P(W^{\Delta})}^{\Sp(W+W^-)} \chi |\det|^s \longrightarrow C^{\infty}(\Sp(W)) \]
and may form the integral
\[  Z(s, \chi,  f_s, v_1, v_2) := \int_{\Sp(W)} f(g,1) \cdot \overline{ \langle gv_1, v_2 \rangle} \, dg,  \quad \text{for $f _s \in I_{P(W^{\Delta})}(s)$ and $v_1, v_2 \in \sigma$.} \]
This is called a local doubling zeta integral. It arises naturally from a global Rankin-Selberg integral by unfolding. 
The following are some basic properties:
\vskip 5pt

\begin{itemize}
\item[(a)]  The integral $Z(s, \chi , f_s, v_1,v_2)$ is absolutely convergent when $Re(s)$ is sufficiently large, and admits a meromorphic continuation to all of $\C$.
This gives a meromorphic family of equivariant linear functionals $Z(s, -)$. 

\item[(b)]  At any $s_0$, the leading term of the Laurent expansion of $Z(s, -)$ at $s= s_0$ defines a nonzero element of
\[  \Hom_{\Sp(W) \times \Sp(W^-)}(I_P(s, \chi) \otimes \overline{\sigma} \otimes \overline{\sigma}^{\vee}, \C). \] 
Recall that this Hom space plays a significant role in the analysis of the doubling seesaw identity.

\item[(c)]   There is a functional equation:
\[   Z(-s, \chi^{-1}, M_{\psi}(s,\chi)(f_s),  v_1, v_2) =   c(s, \sigma, \chi) \cdot  Z(s, \chi, v_1, v_2) \]
where 
\[  M_{\psi}(s ,\chi) : I_P(s,\chi) \longrightarrow I_P(-s, \chi^{-1}) \]
is a (suitably normalized) standard intertwining operator. 
\end{itemize}
Following Tate's thesis, the scalar-valued function $c(s,  \sigma,\chi)$ can be used to give a definition of the 
standard L-factor $L(s, \sigma \times \chi)$ and standard epsilon factor $\epsilon(s, \sigma, \chi,\psi)$ for $G \times \GL_1$ (where $G$ is an arbitrary classical group). This local theory was developed by Lapid and Rallis \cite{LapidRallis2005} and completed by Yamana \cite{Yamana2014} and Kakuhama \cite{Kakuhama2020}.
\vskip 5pt

\subsection{\bf Relevance to epsilon dichotomy}
We can now explain the relevance of the doubling zeta integral to the epsilon dichotomy theorem. When $\sigma$ is tempered, the local doubling zeta integral is convergent at $s = 0$ and
so defines a nonzero element in the 1-dimensional vector  space
\[  \Hom_{\Mp(W) \times \Mp(W^-)}(I_P(0) \otimes \overline{\sigma} \otimes \overline{\sigma}^{\vee}, \C) \]
 in (b) above. Because
 \[  I_{P(W^{\Delta})}(0)  = \bigoplus_{\epsilon= \pm} \Theta_{V^{\epsilon}, W+W^-}(1), \] 
this functional is nonzero on precisely one of the subspaces $\Theta_{V^{\epsilon}, W+W^-}(1) \otimes  \overline{\sigma} \otimes \overline{\sigma}^{\vee}$. To prove the epsilon dichotomy theorem amounts to determining precisely which of the two summands give nonzero contribution. 
Now the normaiized intertwining operator $M_{\psi}(0)$ is an endomorphism of $I_P(0)$ and  acts as $\epsilon$ on the summand $\Theta_{V^{\epsilon}, W+W^-}(1)$. This together with the local functional equation determines precisely which summand has nonzero contribution in terms of the local root number $\epsilon(1/2,\sigma, \psi)$. 
\vskip 5pt

\vskip 5pt

\subsection{\bf The stable range and unitarity}  \label{SS:stable range}
The local doubling zeta integral in facts gives one an explicit construction of the local theta lifting and  appeared in the thesis work of Jianshu Li \cite{Li87}, concurrently and somewhat independently of the work \cite{GPSR87} of Piatetski-Shapiro and Rallis. Let us recall the context of Li's work.
\vskip 5pt

We have seen that if $\Sp(W) \times \O(V)$ is in the stable range, say with $\dim  V \geq 2 \dim W$, then for any $\sigma \in {\rm Irr}(\Sp(W))$, $\theta(\sigma)$ is nonzero on $\O(V)$. 
Using the local doubling zeta integral, J.S. Li showed further that if $\sigma$ is unitary, then $\theta(\sigma)$ is also unitary, so that the local theta correspondence gives an injection
\[  \theta_{\psi}:  \widehat{\Sp(W)} \hookrightarrow \widehat{\O(V)} \]
where $\widehat{G}$ stands for the unitary dual of $G$.  
\vskip 5pt

More precisely, recall that we have the maps due to Rallis:
\[ \Omega_{V,W, \psi} \otimes \overline{\Omega_{V, W,\psi}} \cong  \mathcal{S}(V \otimes W^{\nabla}) \twoheadrightarrow  \Theta_{V, W+W^-}(1)  \hookrightarrow  I_P(s_{V,W}). \]
Here, the first map is given by a partial Fourier transform (often called a Weyl transform). 
If $F$ is nonarchimedean, the last injection is in fact an isomorphism, as the degenerate principal series is irreducible.  
For $\phi_1, \phi_2 \in \Omega_{V,W,\psi}$, with corresponding image $f_{\phi_1,\phi_2} \in  I_P(s_{V,W})$, it follows from definition that
\[  f_{\phi_1,\phi_2}(g) = \langle g \cdot \phi_1, \phi_2 \rangle, \quad \text{for $g \in \Sp(W) \subset \Sp(W + W^-)$.}\]
The local zeta integral is thus given by the convergent integral of the matrix coefficients:
\[  Z(s_{V,W}, \chi_V, f_{\phi_1,\phi_2}, v_1, v_2)  = \int_{\Sp(W)} \langle g \cdot \phi_1, \phi_2 \rangle \cdot \overline{\langle g \cdot v_1, v_2 \rangle} \, dg. \]
This is a $\Sp(W) \times \Sp(W) \times \O(V)^{\Delta}$-equivariant nonzero functional
\[  \Omega_{V,W,\psi} \otimes \overline{\Omega_{V,W,\psi}} \otimes \overline{\sigma} \otimes \sigma \longrightarrow \C \]
which  factors as:
\[  \Omega_{V,W,\psi} \otimes \overline{\Omega_{V,W,\psi}} \otimes \overline{\sigma} \otimes \sigma  \twoheadrightarrow \Theta(\sigma) \otimes \overline{\Theta(\sigma)} \longrightarrow \C. \]
Hence, it gives a nonzero $\O(V)$-invariant Hermitian form on the big theta lift $\Theta(\sigma)$. Li showed \cite{Li87}:
\vskip 5pt

\begin{itemize}
\item the kernel of this Hermitian form is simply  the unique maximal proper submodule of $\Theta(\sigma)$, so that it descends further to give a nondegenerate invariant Hermitian pairing on $\theta(\sigma)$. 

\item this pairing on $\theta(\sigma)$ is positive,  giving an invariant inner product, so that 
$\theta(\sigma)$ is unitary.
\end{itemize}
 As mentioned in \S \ref{SS:big-v-small}, it was shown more recently by Loke-Ma  \cite{LokeMa} that when $F = \R$,  one has $\Theta(\sigma) = \theta(\sigma)$  for unitary $\sigma$ in the stable range. It will be good to know this for nonarchimedean $F$ as well.
\vskip 5pt

Hence, we see that in the stable range, the local doubling zeta integral (or equivalently the integral of matrix coefficients) provides an explicit construction of the local theta lifting and  the invariant inner product on the unitary theta lifts. 
 Li also characterized the image of the theta lifting in the stable range as those unitary representations of low rank, in the sense of Howe.
\vskip 5pt

\subsection{\bf The $L^2$-problem}  \label{SS:L2 problem}
We can now discuss the $L^2$-version of the local theta correspondence. Assuming that $\dim W <  \dim V$ so that $\Sp(W)$ is the smaller group. Then we may write:
\[  \hat{\Omega}_{V,W, \psi} \cong \int_{\widehat{\Sp(W)}} \hat{\sigma} \otimes \Theta(\hat{\sigma}) \, d\mu(\hat{\sigma}), \]
for some measure $d\mu$ on the unitary dual of $\Sp(W)$. The goal is to explicate the above isomorphism, the measure $d\mu$ and the unitary representation $\Theta(\hat{\sigma})$ for $d\mu$-almost all $\sigma$. 

\vskip 5pt

In \cite{Sake2017}, Sakellaridis showed:
\vskip 5pt

\begin{thm}  \label{T:L2 problem}
One has a spectral decomposition of the invariant inner product on $\hat{\Omega}_{V,W,\psi}$:
\[  \langle \phi_1, \phi_2 \rangle = \int_{\widehat{\Sp(W)}}  J_{\hat{\sigma}} (\phi_1, \phi_2) \,  d\mu_{\Sp(W)}(\hat{\sigma}) \]
where $d\mu_{\Sp(W)}$ is the (Harish-Chandra) Plancherel measure   on $\widehat{\Sp(W)}$ (relative to a fixed Haar measure $dg$ on $\Sp(W)$) and
\begin{align}
 J_{\hat{\sigma}}(\phi_1 , \phi_2) &=   \sum_{v \in \mathcal{B}(\hat{\sigma})}  \int_{\Sp(W)} \langle g \cdot \phi_1, \phi_2 \rangle \cdot \overline{\langle g \cdot v_1, v_2 \rangle} \, dg. \notag \\
 &= \sum_{v \in \mathcal{B}(\hat{\sigma})} Z(s_{V,W}, \chi_V, \phi_1, \phi_2, v, v) \notag 
 \end{align}
where $\mathcal{B}(\hat{\sigma})$ stands for an orthonormal basis of $\hat{\sigma}$. 
\end{thm} 
Hence, we see that:
\vskip 5pt

\begin{itemize}
\item the spectral support of $\hat{\Omega}_{V,W,\psi}$, as a representation of $\Sp(W)$, is contained in the tempered dual of $\Sp(W)$. 
\item the spectral decomposition of the invariant inner product is basically given by the local doubling zeta integral, or equivalently by the integral of matrix coefficients considered by J.S. Li.
\end{itemize}
The proof of the theorem amounts to the observation that a matrix coefficient $f_{\phi_1, \phi_2} = \langle  - \cdot \phi_1, \phi_2 \rangle$ of the Weil representation $\hat{\Omega}_{V,W,\psi}$, when regarded as a function   on $\Sp(W)$, belongs to the Harish-Chandra Schwarz space. One can then perform a spectral expansion of $f_{\phi_1, \phi_2}$  using the Harish-Chandra Plancherel theorem. Note that $J_{\sigma}(\phi_1,\phi_2)$ 
can also be interpreted as (the complex conjugate of) the  trace of the  convolution operator induced by $f_{\phi_1, \phi_2}$ on the tempered smooth representation $\sigma = \hat{\sigma}^{\infty}$. Hence, we can also express the above theorem as:
\[  \langle \phi_1, \phi_2 \rangle = \int_{\widehat{\Sp(W)}} \overline{{\rm Tr}(\sigma(f_{\phi_1, \phi_2}))} \, \,  d\mu_{\Sp(W)}(\hat{\sigma}). \]
\vskip 10pt

To summarize, we have highlighted the important role played by the doubling seesaw and the doubling zeta integral in many aspects of the theory of local theta correspondence.

\vskip 10pt
\section{\bf Global Results}
Let us now consider the global setting where we work over a number field $k$ with adele ring $\A$. We shall highlight some results on the global theta correspondence.
\vskip 5pt

\subsection{\bf Global theta liftings}  \label{SS:global theta2}
 We have introduced the``formation of theta series":
\[  \begin{CD}
\Omega_{V,W, \psi} = \mathcal{S}(\mathbb{Y}_{\A})  @>\theta>>   \mathcal{A}(\Mp(\mathbb{W})) @>i_{V,W,\psi}^*>>  \mathcal{C}([\O(V) \times \Sp(W)]) \end{CD} \]
 and the global theta lifting for a cuspidal representation
   $\sigma \subset \mathcal{A}_{cusp}(\Sp(W))$. Namely,  for $\phi \in \Omega_{V,W,\psi}$ and $f \in \sigma$, we have
\[  \theta(\phi,f)(h) := \int_{[\Sp(W)]} \theta(\phi) (gh) \cdot \overline{f(g)} \, dg, \]
and
\[  \Theta(\sigma) := \langle \theta(\phi,f): \phi \in \Omega_{V,W,\psi}, \, f \in \sigma \rangle \subset \mathcal{A}(\O(V)). \]
 We have also raised the following questions:
\begin{itemize}
\item Decide if $\Theta(\sigma)$ is nonzero.
\item Is  $\Theta(\sigma) \subset \mathcal{A}_{cusp}(\O(V))$? Is $\Theta(\sigma) \subset \mathcal{A}_2(\O(V))$?
\item What is the relation between the global theta lift and the abstract theta lift 
\[  \theta^{abs}(\sigma) := \bigotimes_v \theta(\sigma_v) ?\]
\end{itemize}
We will discuss these questions in reverse order, highlighting the complete parallel in the treatment of these global problems with that of the local ones. 
\vskip 5pt

\subsection{\bf Compatibility with local theta lifts}  \label{SS:compatibility}
For the third question above, one has \cite{Rallis84a, Gan08}:
 \begin{prop}
 Suppose that $\Theta(\sigma)$ is nonzero and is contained in the space $\mathcal{A}_2(\O(V))$ of square-integrable automorphic forms on $\O(V)$. 
Then $\Theta(\sigma) \cong \theta^{abs}(\sigma)$.
\end{prop}
This proposition is a direct consequence of the (local) Howe duality theorem. Under the hypotheses of the proposition, the determination of the isomorphism class of $\Theta(\sigma)$ is reduced to a local problem.
\vskip 10pt

\subsection{\bf Rallis tower}  \label{SS: global rallis tower}
As in the local case, for the questions of cuspidality and nonvanishing, 
it is useful to consider a Rallis tower of theta lifitings, corresponding to a Witt tower $V_r = V_0 \oplus \mathbb{H}^r$ of quadratic spaces, with $V_0$ anisotropic. 
One has the analogous results \cite{Rallis84a, Rallis84b}:
\vskip 5pt

\begin{prop}  \label{P:global tower}
Let $\sigma$ be a cuspidal automorphic representation of $\Sp(W)$, and consider its global theta lift $\Theta_{V, W_r,\psi}(\sigma)$ on $\O(V_r)$. Then one has:
\vskip 5pt

(i) There is a smallest  $r_0 = r_0(\sigma) \leq \dim W$ such that $\Theta_{V_{r_0}, W, \psi}(\sigma) \ne 0$.  Moreover,  $\Theta_{V_{r_0} W, \psi}(\sigma) $ is contained in the space of cusp forms.

\vskip 5pt

(ii) For all $r >  r_0$, $\Theta_{V_r, W, \psi}(\pi)$ is nonzero and noncuspidal.

\vskip 5pt

(iii) For all $r \geq \dim W$, $\Theta_{V_{r_0}, W, \psi}(\sigma)  \subset \mathcal{A}_2(\O(V_r))$. 
\end{prop}
As in the local case, we call $r_0= r_0(\sigma)$ the first occurrence of $\pi$ in the relevant Witt tower, and the range where $r \geq \dim W$ is called the stable range.  
\vskip 10pt

 \subsection{\bf Rallis inner product}  \label{SS:inner product SW}
 To detect the nonvanishing of $\Theta(\sigma) \subset \mathcal{A}(\O(V))$, one may assume (in view of Proposition \ref{P:global tower}) that the theta lift of $\sigma$ to the lower step of the Rallis tower is zero, so that $\Theta(\sigma) \subset \mathcal{A}_{cusp}(\O(V))$. Then one may check the nonvanishing or not of $\theta(\phi,f) \in \Theta(\sigma)$ by computing 
 its Petersson inner product. This computation was first considered by Rallis in \cite{Rallis87} and   is based on the doubling seesaw:
 \[
    \xymatrix{
    \Sp(W \oplus W^-)\ar@{-}[d]\ar@{-}[dr]& \O(V)\times \O(V)\ar@{-}[d]\\
    \Sp(W)\times \Sp(W^-)\ar@{-}[ur]& \O(V)^{\Delta}
    }
\]
 and is the global analog of the local seesaw identity.
 \vskip 5pt

For $\phi_1, \phi_2 \in \Omega_{V,W,\psi}$ and $f_1, f_2 \in \sigma$,  one has:
\begin{align}  \label{A:inner}
 & \langle \theta(\phi_1, f_1), \theta(\phi_2,f_2) \rangle  \notag \\
 = & \int_{[\O(V)]} \left( \int_{ [\Sp(W)]} \theta(\phi_1)(g_1,h)
   \cdot \overline{f_1(g_1)} \, dg_1 \right) \cdot 
\left(  \int_{ [\Sp(W)]} 
\overline{ \theta(\phi_2)(g_2,h)} \cdot f_2(g_2) \, dg_2 \right) \, dh \notag \\
= & \int_{ [\Sp(W) \times  \Sp(W)]} \left( \int_{ [\O(V)]}
  \theta(\phi_1)(g_1, h) \cdot \overline{ \theta(\phi_2) (g_2,h)} \,
  dh \right) \cdot  
  \overline{f_1(g_1)} \cdot f_2(g_2) \, dg_1 \, dg_2 
\end{align}
where in the last equality, we have formally exchanged the integrals. 
This inner integral $I(\phi_1 \otimes \phi_2)$  (if it converges) can be interpreted as the global theta lift of the constant function $1$ of $\O(V)^{\Delta}$ to $\Sp(W \oplus W^-)$.
 The reader will recognise that this use of the Fubini theorem is the global manifestation of the local seesaw identity.
 
 \vskip 5pt
 
 \subsection{\bf Siegel-Weil formula}
To proceed further, one would like to give a different interpretation of the inner integral 
\[  I(\phi_1\otimes \phi_2) = \int_{ [\O(V)]}
  \theta(\phi_1)(g_1, h) \cdot \overline{ \theta(\phi_2) (g_2,h)} \,dh. \]
  The reader will recognise this as the global analog of the local determination of $\Theta(1)$. 
  The understanding of this inner integral is the content of the so-called {\em Siegel-Weil formula}: it identifies the inner integral with an Eisenstein series.
 More precisely, recall that one has a sequence of maps
 \[  \Omega_{V,W, \psi} \otimes \Omega_{V, W^-, \psi} \cong \Omega_{V, W+W^-,\psi} \twoheadrightarrow \Theta_{V, W+W^-,\psi}(1) \hookrightarrow I_{P(W^{\Delta})}(s_{V,W}) \]
 For $\phi_1, \phi_2 \in  \Omega_{V,W,\psi}$, let 
 \[  f_{\phi_1, \phi_2} \in  I_{P(W^{\Delta})}(s_{V,W}, \chi_V)  \]
 be its image under the above sequence of maps. Then to a first degree of approximation, the Siegel-Weil formula says:
 \[   I(\phi_1 \otimes \phi_2) = E(s_{V,W}, f_{\phi_1, \phi_2}), \]
 where the RHS denotes the Eisenstein series associated to $f_{\phi_1, \phi_2}$. As mentioned, this is just a first approximation: the theta integral on the LHS may need regularization as it may diverge and the Eisenstein series on the RHS may have pole at $s =s_{V,W}$ and so the RHS needs to be suitably interpreted.  Deriving a precise Siegel-Weil formula is a subject that takes close to 40 years, beginning with the pioneering work of Weil \cite{Weil65}, followed by the breakthrough papers of Kudla and Rallis \cite{KudlaRallis88a, KudlaRallis88b, KudlaRallis94}, refined  in the work of Ikeda \cite{Ikeda96}, Ichino \cite{Ichino2001, Ichino2004, Ichino2007} and Yamana \cite{Yamana2011, Yamana2013a, Yamana2013b},  before culminating in the paper of Gan-Qiu-Takeda \cite{GQT2014}. 
 
\vskip 5pt

\subsection{\bf Global doubling zeta integral}
In any case, continuing with our informal discussion, we apply the Siegel-Weil formula  to obtain 
\[  \langle \theta(\phi_1, f_1), \theta(\phi_2,f_2) \rangle  = \int_{ [\Sp(W) \times  \Sp(W)]}  E(s_{V,W}, f_{\phi_1, \phi_2},g) \cdot  \overline{f_1(g_1)} \cdot f_2(g_2) \, dg_1 \, dg_2. \]
The integral
\[  Z(s, \chi_V, f_s,  f_1, f_2) = \int_{ [\Sp(W) \times  \Sp(W)]}  E(s, f_s,g) \cdot  \overline{f_1(g_1)} \cdot f_2(g_2) \, dg_1 \, dg_2. \] 
is the global doubling zeta integral. When $Re(s)$ is sufficiently large, it unfolds to give an Euler product of the local doubling zeta integral encountered earlier:
\[  Z(s, \chi_V, f_s,  f_1, f_2) = \prod_v Z_v(s, \chi_{V,v}, f_{s,v}, f_{1,v}, f_{2,v}) \]
As mentioned before, this Rankin-Selberg integral was discovered by Piatetski-Shapiro and Rallis \cite{GPSR87}, who showed that it represents the standard L-function of $\Sp(W) \times \GL_1$:
\[  Z(s, \chi_V, f_s,  f_1, f_2) =  \frac{L^S(s+ \frac{1}{2}, \sigma \times \chi_V)}{\alpha^S(s)} \cdot \prod_{v \in S}  Z_v(s, \chi_{V,v}, f_{s,v}, f_{1,v}, f_{2,v}) \]
where $S$ is a sufficiently large finite set of places outside of which every data involved is unramified and $\alpha^S(s)$ is a product of abelian Hecke L-functions.
The local theory of the doubling zeta integral was completed by Lapid-Rallis \cite{LapidRallis2005}, Yamana \cite{Yamana2014} and Kakuhama \cite{Kakuhama2020}, and as we saw earlier, gave rise to the definition of the local L-factors and epsilon factors at all places. 
\vskip 5pt

As a summary of the above discussion, we obtain \cite{GQT2014}:
\vskip 5pt

\begin{thm} \label{T:rallis inner product}
Let $\sigma \subset \mathcal{A}_{cusp}(\Sp(W))$ be such that $\Theta(\sigma) \subset \mathcal{A}_{cusp}(\O(V))$. Assume (for simplicity) that $\dim W < \dim V$, so that 
\[ s_{V,W} = \frac{\dim V - \dim W -1}{2} \geq 0.\]
\vskip 5pt
 
(i) (Rallis inner product formula)  For  $\phi_1, \phi_2 \in \Omega_{V,W,\psi}$ and $f_1, f_2 \in \sigma$, one has
\[  \langle \theta(\phi_1, f_1), \theta(\phi_2,f_2) \rangle  =  \frac{L(s_{V,W}+ \frac{1}{2}, \sigma \times \chi_V)}{\alpha(s_{V,W})} \cdot \prod_v  Z^*_v(s_{V,W}, \chi_{V,v}, f_{\phi_1, \phi_2}, f_{1,v}, f_{2,v})\]
where $Z_v^*$ refers to a normalized local doubling zeta integral and $\alpha(s)$ is a product of abelian Hecke L-functions.
\vskip 5pt

(ii) (local-global criterion for nonvanishing)  The global theta lift $\Theta(\sigma)$ is nonzero if and only if the following two conditions hold:
\vskip 5pt

\begin{itemize}
\item[(a)] $L(s_{V,W}+ \frac{1}{2}, \sigma \times \chi_V) \ne 0$;
\vskip 5pt

\item[(b)]  for all $v$, $Z^*_v(s_{V,W}, -)$ is a nonzero functional on $\Omega_{V_v, W_v, \psi_v} \otimes \Omega_{V_v, W^-_v, \psi_v} \otimes \overline{\sigma_v} \otimes \sigma_v$.
\end{itemize}
Moreover, for each finite or complex $v$, the condition (b) above is equivalent to the nonvanishing of the local theta lift $\theta(\sigma_v)$. 
\end{thm}
Naturally, one expects the last equivalence to hold at real places $v$ as well; this will complete the local-global vanishing criterion. This is known in some cases over $\R$ and it would be good to complete the story here.
\vskip 5pt

 We stress again that the global argument in this subsection is the analog of the local argument discussed in \S \ref{SS:doubling seesaw} and \S \ref{SS:proof of Howe}.
\vskip 10pt

\section{\bf Contact with the Langlands Program}

We have now spent quite a number of pages discussing various fundamental results within the theory of theta correspondence, highlighting especially the roles of the doubling seesaw and  the doubling zeta integrals, in both the local and global settings. All these results were developed to address the basic questions of (Nonvanishing) and (Identity) mentioned in \S \ref{SS:further questions}. While these results help clarify the situation, they do not fully address these questions. 
It is at this point that the Langlands program enters the picture.
\vskip 5pt
 Indeed,  the local Langlands correspondence  provides a classification of irreducible representations of real or p-adic reductive groups,  akin to the highest weight theory for connected compact Lie groups. It thus provides a language for us to say when $\theta(\pi)$ is nonzero and what $\theta(\pi)$ is.   
 From this point of view,  the Langlands program is providing a service to the theory of theta correspondence.
 
 \vskip 5pt

\subsection{\bf Theta lifts of unramified representations} \label{SS:unramified reps}
 Let us illustrate this in the case of unramified representations. Let ${\rm Irr}_{K_W}(\Sp(W))$ be the set of $K_W$-unramified representations of $\Sp(W)$, where $K_W$ is a hyperspecial maximal compact subgroup of $\Sp(W)$ stabilizing a self-dual lattice $\Lambda_W$ in $W$. If $\mathcal{H}(\Sp(W), K_W)$ is the associated spherical Hecke algebra (which is commutative), then the Satake isomorphism gives rise to a bijection
 \[  {\rm Irr}_{K_W}(\Sp(W)) \leftrightarrow {\rm Spec}(\mathcal{H}(\Sp(W), K_W)) =   \{ \text{semisimple conjugacy classes in $\Sp(W)^{\vee}$\},} \]
 attaching to an  unramified representation $\sigma$ its Satake parameter $s_{\sigma}$. 
 One might say that this bijection is one of the impetus that initiated the Langlands program. Likewise, assuming $V$ to be a split quadratic space for simplicity, one has a bijection
 \[  {\rm Irr}_{ur}(\O(V)) \leftrightarrow {\rm Spec}(\mathcal{H}(\O(V), K_V)) = \{ \text{semisimple conjugacy classes in $\O(V)^{\vee}$} \} \]
 for $K_V$-unramified representations of $\O(V)$,  where $K_V$ is a hyperspecial maximal compact subgroup stabilizing a self-dual lattice $\Lambda_V$ in $V$. 
 \vskip 5pt
 
  Let $\psi : F \rightarrow \C^{\times}$ have conductor $\mathcal{O}_F$ and consider the associated Weil representation $\Omega_{V,W,\psi}$. If $K$ is the stabilzer of the lattice $\Lambda_V \otimes \Lambda_W$ in  $\Sp(V \otimes W)$, then 
  \[  \Omega_{V,W, \psi}^K = \C \cdot \phi_0 \]
  is 1-dimensional.  Now one has the following classic result of Howe (see \cite[Pg. 103 and Pg 105]{MVW87}) and Rallis \cite[\S 6]{Rallis82}:
 \vskip 5pt
 
 \begin{thm}  \label{T:unramified}
  In the above context, 
  \[  \Omega_{V,W,\psi}^{K_W} = \text{$C^{\infty}_c(O(V))$-span of $\phi_0$}  \] 
 and
 \[  \Omega_{V,W,\psi}^{K_V} = \text{$C^{\infty}_c(\Sp(W))$-span of $\phi_0$}. \]
 Further, one has:
 \vskip 5pt
 
 \begin{itemize}
 \item[(a)]  If $\pi \otimes \sigma \in {\rm Irr}(\O(V) \times \Sp(W)$ is a quotient of $\Omega_{V,W,\psi}$, then $\pi$ is $K_V$-unramified if and only if $\sigma$ is $K_W$-unramified.
Hence the theta correspondence gives a map
\[  {\rm Irr}_{K_W}(\Sp(W)) \longrightarrow {\rm Irr}_{K_V}(\O(V))  \sqcup \{0 \}. \]
 \vskip 5pt
 
 \item[(b)]   Assume without loss of generality that $\dim W \leq  \dim V -1$. There is a surjective  algebra homomorphism 
 \[  f:   \mathcal{H}(\O(V), K_V)   \longrightarrow \mathcal{H}(\Sp(W), K_W) \]
 such that  the map in (a) above is given by the induced map $f^*$ on the maximal spectra. In particular, the theta lift of an unramified $\sigma \in {\rm Irr}(\Sp(W))$ is nonzero on $\O(V)$.
  
 \vskip 5pt
 \item  Suppose $\dim W = 2n$ and $\dim V = 2m > 2n$, so that  $\Sp(W)^{\vee}  = \SO_{2n+1}(\C)$ and $\O(V)^{\vee} = \O_{2m}(\C)$. 
  In terms of Satake parameter, the map 
  \[  f^* :   \SO_{2n+1}(\C)_{ss} /_{\text{conj.}}  \longrightarrow \O_{2m}(\C)_{ss}/_{\text{conj.}}. \]
    is given explicitly by:
 \[  f^*(s) =  s \oplus  \Sym^{2m-2n-2} \left( \begin{array}{cc} 
 q^{1/2} & \\
 & q^{-1/2} \end{array} \right)   \in \SO_{2n+1}(\C) \times \SO_{2m-2n-1}(\C)  \subset \O_{2m}(\C), \]
 where 
   $q$ is the size of the residue field of $F$. 
  \end{itemize}
 \end{thm}
 
 \vskip 5pt
 
 This theorem illustrates how having a classification of the unramified representations (i.e. the unramified local Langlands correspondence) greatly facilitates the explicit description of the local theta correspondence.  Thanks to advances in the Langlands program for classical groups in recent years, both the questions of (Nonvanishing) and (Identity) have now rather complete answers.   We refer the reader to \cite{AtobeGan2017, BakicHanzer2021} for the precise statements. 
 Instead of discussing these, we will next discuss how the theory of theta correspondence returns the favour by contributing to the Langlands program.

\vskip 10pt

\subsection{\bf Theta lifts and A-parameters}

 In the global setting, it turns out that the global theta correspondence is best described using the language of A-parameters. (We assume familiarity with the notion of A-parameters here and refer the reader to Kaletha's lectures for further details of this and the Arthur conjecture).
For example, we saw in Proposition \ref{P:global tower} that, in the stable range,  the global theta lift of a cuspidal $\sigma$ is always nonzero and contained in $\mathcal{A}_2(\O(V))$. 
 On the other hand, one knows that the spectral decomposition of $\mathcal{A}_2(G)$ is governed by Arthur's conjecture, which is a classification of the near equivalence classes by elliptic A-parameters of $G$. Hence,    it is natural to ask how the A-parameters of $\Theta(\sigma)$ and $\sigma$ are related whenever $\Theta(\sigma)$ is square-integrable. 
 \vskip 5pt
 
 If we view A-parameters as parametrizing near equivalence classes, answering this question involves understanding the local theta lifts of unramified representations , which is precisely what we have done in Theorem \ref{T:unramified} above.    Indeed, Theorem \ref{T:unramified} motivates: 
 \vskip 5pt
 
 \begin{con}[Adam's conjecuture \cite{Adams89}]  \label{conj:adams}
 If $\sigma \subset \mathcal{A}_{cusp}(\Sp(W))$ has A-parameter $\Psi$ (thought of as an $\dim W +1$-dimensional representation of $L_F \times \SL_2(\C)$), then the global theta lift of $\sigma$ (which we assume is a summand in $\mathcal{A}_2(\O(V))$) has A-parameter
 \[    \Psi  \oplus  S_{\dim V - \dim W -1}. \]
 \end{con}
 
Given Arthur's conjecture, this is largely a local unramified issue.  What is subtle about Adam's conjecture \cite{Adams89} is its local analog (which we won't discuss here).
The investigation of the Adams conjecture in the local context was initiated by Moeglin \cite{Moeglin2011b} and recently completed by Baki\'c-Hanzer \cite{BakicHanzer2022}.
  \vskip 5pt

\subsection{\bf Applications}
Taking Adams' conjecture as a guide, we can appreciate how  theta correspondence is often used to produce interesting examples of square-integrable automoprhic forms on a classical group $G$ as liftings from cusp forms on smaller groups $H$.  Here are 3 well-known examples:
\vskip 5pt
\begin{itemize}
\item  the theta lifting from $\O_2$ to $\Sp_4$ was used by Howe and Piatetski-Shapiro \cite{HPS79} to produce the first counterexamples to the naive generalized Ramanujan conjecture in the Corvallis proceedings.
\vskip 5pt

\item the theta lifting for $\Mp_2 \times \O_3$ was used by Waldspurger \cite{Wald80} to give an automorphic rendition of the Shimura correspondence, leading to a beautiful classification of cuspidal representations of $\Mp_2$.

\vskip 5pt

\item the theta lifting from $\Mp_2$ to $\O_5$ was used by Piatetski-Shapiro \cite{PS83} to construct and understand the Saito-Kurolawa automorphic representations of $\PGSp_4 \cong \SO_5$.
\end{itemize}
\vskip 5pt

\noindent  In view of Adam's conjecture, one would not expect such construction to yield the whole automorphic discrete spectrum of the larger group $G$. Indeed, Adams conjecture predicts that one would only obtain A-parameters on $G$  of the following form:
\[  \text{(A-parameters of smaller group)} \oplus (\chi \otimes S_r ) \]
for some 1-dimensional character $\chi$. Such A-parameters are of course far from being the most general. This is of course not surprising:  there is no reason to expect that the automorphic discrete spectrum of a group can be totally understood in terms of the spectra  of smaller groups.  
\vskip 5pt

\subsection{\bf Work of Chen-Zou}
In the rest of this section, we will describe a series of recent work by Rui Chen and Jialiang Zou \cite{ChenZou2021}, which (somewhat surprisingly) establishes Arthur's conjecture (at least the tempered part) for non-quasi-split classical groups, using the theta correspondence to transport one's knowledge from the quasi-split case established by Arthur and Mok. To be more precise, the approach  was first conceived in \cite{GanIchino2018} to show  partially the Arthur conjecture for metaplectic groups $\Mp(W)$, before it was adapted to the case of non-quasi-split classical groups in the thesis work of Chen and Zou.
\vskip 5pt

\subsection{\bf The idea}
We shall illustrate with the case of even orthogonal groups.
The main idea is simple: instead of lifting from smaller groups to a fixed (non-quasi-split) orthogonal group $\O(V)$, why don't we construct $\mathcal{A}_{2}(\O(V))$ by theta lifting from a larger symplectic group $\Sp(W)$ (which is split)? More precisely, if you would like to construct and understand  $\mathcal{A}_{\psi}(\O(V))$  (the near equivalence class in $\mathcal{A}_2(\O(V))$ determined by $\psi$) for a particular A-parameter $\psi$, you can do so as follows:
\vskip 5pt

\begin{itemize}
\item take a much larger $\Sp(W)$, so that $\O(V) \times \Sp(W)$ is in the stable range, and consider the A-parameter 
\[  \psi' = \psi \oplus  S_r \]
of $\Sp(W)$ (as suggested by Adam's conjecture). By the work of Arthur \cite{Arthur2012}, one already ``knows" the submodule $\mathcal{A}_{\psi'}(\Sp(W))$ in terms of local and global A-packets with a multiplicity formula. 

\vskip 5pt

\item use global theta lifting to transport $\mathcal{A}_{\psi'}(\Sp(W))$ back to $\O(V)$; one should get the submodule $\mathcal{A}_{\psi}$.
\end{itemize}

\vskip 5pt

\subsection{\bf Some Issues}
There will however be some potential issues that one can readily point out:
\vskip 5pt

\begin{itemize}
\item Since the A-parameter $\psi'$ of $\Sp(W)$ is nontempered, the understanding of the associated A-packet coming from Arthur's work is far from complete. For example, Arthur does not know from his approach how many representations are in the local A-packets. If one hopes to construct $\mathcal{A}_{\psi}$ from $\mathcal{A}_{\psi'}$, this is  undesirable:  one cannot hope to have a good understanding of $\mathcal{A}_{\psi}$ without first understanding the input $\mathcal{A}_{\psi'}$. 
\vskip 5pt

Thankfully, this issue is somewhat alleviated by the independent work of Moeglin \cite{Moeglin2009, Moeglin2011a}, which was clarified and extended by Bin Xu \cite{Xu2017, Xu2021}, Hiraku Atobe \cite{Atobe2022a, Atobe2022b} and Heseltine-Liu-Lo \cite{HazeltineLiuLo2022}.  As a result, one has an independent and explicit construction of local A-packets using Jacquet module techniques. In particular,   one has a rather good understanding of $\mathcal{A}_{\psi'}$. 
\vskip 10pt

\item By the tower property of theta correspondence, one does not expect many constituents of $\mathcal{A}_{\psi'}$ to be cuspidal; indeed, they should be noncuspidal when $\Sp(W)$ is sufficiently large. Given this, how can one perform the global theta lifting on them? Recall that we need the integral defining the global theta lifting to converge, and there may be no reason for this convergence when the input is noncuspidal. 
 \end{itemize}
\noindent In the next subsection, we shall explain the novel idea that allows one to bypass this last difficulty.

\vskip 5pt

\subsection{\bf J.S. Li's Inequalities}
In the stable range, it turns out that one can study the global  theta lifting from the $L^2$-point of view, exploiting the fact that the representations on the larger group has low rank. 
This innovative approach was carried out by J.S. Li in his 1997 paper \cite{Li97}. Let us recall his results here.
\vskip 5pt

\begin{theorem}  \label{T:Lis}
Assume that $\O(V) \times \Sp(W)$ is in the stable range with $\O(V)$ the smaller group.
Suppose that one has a direct integral decomposition
\[  L^2([\O(V)]) = \int_{\widehat{\O(V)(\A)}}   m_{\pi} \cdot \pi\, d\mu(\pi) \]
for some measure $d\mu$ on the unitary dual $\widehat{\O(V)(\A)}$ of $\O(V)(\A)$, and some multiplicity function $m(-)$.
Then $L^2([\Sp(W)])$ contains a submodule isomorphic to the direct integral
\[   \int_{\widehat{\O(V)(\A)}}   m_{\pi} \cdot \theta^{abs}(\pi)\, d\mu(\pi) \]
where $\theta^{abs}(\pi) = \otimes_v \theta(\pi_v)$ is the (abstract) theta lift  of $\pi$ (which is nonzero since we are in the stable range). 
\end{theorem}
\vskip 5pt
It is worth noting that the proof of this theorem does not involve the usual integral defining the global theta lifting, but involves harmonic analysis of the part of the $L^2$ spectrum of $\Sp(W)$ involving low rank representations. The properties of low rank representations thus play an indispensable role here.
\vskip 5pt

As a corollary of this theorem and other results, we have the following inequalities:
\vskip 5pt

\begin{cor}  \label{C:Lis}
For $\pi \in {\rm Irr}(\O(V)(\A))$, let
\[  m(\pi) = \dim \Hom_{\O(V)}(\pi, \mathcal{A}(\O(V))) \quad \text{and} \quad m_2(\pi) =  \dim \Hom_{\O(V)}(\pi, \mathcal{A}_2(\O(V))). \]
Then one has
\[  m_2(\pi) \leq m_2(\theta^{abs}(\pi)) \leq m(\theta^{abs}(\pi)) \leq m(\pi). \]
\end{cor}
Here, the first inequality is a consequence of the theorem, the second inequality is obvious whereas the third  follows by other considerations (of Fourier coefficients).

\vskip 5pt

\subsection{\bf Assigning A-parameters}
Let us illustrate how the corollary above is useful in attaching A-parameters to near equivalence classes. Suppose that $\mathcal{C} \subset \mathcal{A}_2(\O(V))$ is a nonzero near equivalence class, say
\[  \mathcal{C} \cong \bigoplus_{i \in I} m(\pi_i) \pi_i. \]
By the corollary,  we see that $\mathcal{A}_2(\Sp(W))$ contains as a submodule
\[  \theta^{abs}(\mathcal{C}) := \bigoplus_{i \in I} m(\pi_i) \cdot \theta^{abs}(\pi_i). \]
Now all the summands in $\theta^{abs}(\mathcal{C})$ are nonzero and nearly equivalent to each other. Since $\Sp(W)$ is split, one knows by Arthur that $\theta^{abs}(\mathcal{C})$ is associated to an elliptic A-parameter $\psi'$.  By using our knowledge of the unramified theta correspondence (as given in Theorem \ref{T:unramified}) and poles of standard L-functions (coming from the doubling zeta integrals), one can show that $\psi'$ has the desired form
\[ \psi' =  \psi \oplus  S_r \]
for some elliptic A-parameter $\psi$ for $\O(V)$. One then defines $\psi$ to be the A-parameter of the near equivalence class  $\mathcal{C}$. 
\vskip 5pt

In this way, Chen and Zou \cite{ChenZou2021} showed:
\vskip 5pt

\begin{theorem}  \label{T:CZ}
For any (non-quasi-split) even orthogonal group $\O(V)$, one has a decomposition
\[  \mathcal{A}_2(\O(V)) = \bigoplus_{\psi}  \mathcal{A}_{\psi} \]
where $\psi$ runs over the elliptic A-parameters of $\O(V)$ and $\mathcal{A}_{\psi}$ is the associated near equivalence class.
\end{theorem}
\vskip 5pt
\noindent Another way to formulate this theorem is that it gives the weak Langlands functorial lifting from $\O(V)$ to $\GL(V)$, with image given by the expected description.
\vskip 5pt

\subsection{\bf Does equality hold?}
It remains then to understand each near equivalence class $\mathcal{A}_{\psi}$, and in particular to describe it in the language of Arthur's conjecture. 
In some sense, one would like to transport the structure of $\mathcal{A}_{\psi'}$ back to $\mathcal{A}_{\psi}$. However, since we only have inequalities in the above corollary, it means that one only has
\[  \theta^{abs}(\mathcal{A}_{\psi}) \subset \mathcal{A}_{\psi'}. \] 
In other words, in transferring from $\Sp(W)$ back to $\O(V)$, we might lose some information. For example, it is possible that $\mathcal{A}_{\psi} = 0$
but $\mathcal{A}_{\psi'} \ne 0$. 
\vskip 5pt

Thus, we see that it is important to know when the equality $m_2(\pi) = m_2(\theta^{abs}(\pi))$ holds. For this, one has:
\vskip 5pt

\begin{prop}
Equality holds in J.S. Li's inequalities in Corollary \ref{C:Lis}  in the following cases:
\vskip 5pt

\begin{itemize}
\item $\psi$ is a tempered A-parameter;
\vskip 5pt

\item The Witt index of $V$ is $0$ or $1$. 
\end{itemize}
\end{prop}

In these cases, what was shown is that
\[  m_2(\pi) = m(\pi), \]
so that equalities hold throughout  Li's inequalities in Corollary \ref{C:Lis}.
For these cases, it is then reasonable to expect that the structure of $\mathcal{A}_{\psi}$ can be faithfully inherited from that of $\mathcal{A}_{\psi'}$.

\vskip 5pt

\subsection{\bf Results of Chen-Zou}
In their paper \cite{ChenZou2021}, Chen and Zou showed the following theorem:
\vskip 5pt

\begin{theorem}
(i)  In the context of Theorem  \ref{T:CZ}, the submodule $\mathcal{A}_{\psi}$ can be described as in Arthur's conjecture for any tempered $\psi$.
\vskip 5pt

(ii) Arthur's conjecture holds for $\O(V)$ if $\O(V)$ has $k$-rank $\leq 1$. 
\end{theorem}
\vskip 5pt

To go beyond this, it seems one needs to answer the following question:
\vskip 5pt

\noindent{\bf Question:} Is it the case that one always has: 
\[  m_2(\pi) = m_2(\theta^{abs}(\pi)) = m(\theta^{abs}(\pi)) = m(\pi)? \]
 \vskip 10pt

\section{\bf Exceptional Theta Correspondence}
We have so far focused on the classical setting of theta correspondence, involving dual pairs in symplectic groups. 
Since the 1990's, people have been looking to extend this theory to a more general setting. One can mention, for example, the work of Brylinski-Kostant \cite{BryKostant94a,  BryKostant94b, BryKostant95}, Kazhdan-Savin \cite{Kazhdan90, KazhdanSavin90, Savin94}, Huang-Pandzic-Savin \cite{HPS96}, Ginzburg-Rallis-Soudry \cite{GRS97a,GRS97b} and J.S. Li \cite{Li99}.  As such extensions to classical groups can frequently be  related to the case of symplectic groups, the main novelty is really in the setting of exceptional group, hence the name ``exceptional theta correspondence". 
We give a very brief summary of the state of affairs here, without attempting to be exhaustive.
\vskip 5pt

\subsection{\bf Exceptional dual pairs}
For any reductive group $\mathcal{E}$, one can seek to classify the dual pairs $G \times H \rightarrow \mathcal{E}$. 
On the level of Lie algebra over $\C$, such a classification has been carried out by Rubenthaler \cite{Ruben94}. On the level of groups, one thus needs to figure out the relevant isogeny type and there is also the question of rational forms of these dual pairs over local and global fields.

\vskip 5pt

On the other hand, one can proceed  by a geometric algebra approach. 
Recall, for example, that classical dual pairs were constructed as isometry groups of Hermitian or skew-Hermitian modules (such as symmetric bilinear forms and symplectic forms). 
In a similar vein, dual pairs in exceptional groups can be constructed as automorphism groups of certain algebraic structures. These algebraic structures are those which intervene in the construction of exceptional groups. We highlight two examples here:
\vskip 5pt

\begin{itemize}
\item Composition algebras $C$: these are the base field $F$, a quadratic \'etale $F$-algebra $K$, a quaternion $F$-algebra $B$ or an octonion algebra $\mathbb{O}$.  So $\dim C = 1$, $2$ $4$ or $8$.
\vskip 5pt
\item Freudenthal-Jordan algebras $J$: these are \'etale cubic $F$-algebras $E$ or  the algebra of Hermitian matrices with entries in a composition algebra $C$ and their rational forms. So $\dim J = 3$, $6$, $9$, $15$ or $27$. 
\end{itemize}

 Given a composition $F$-algebra $C$ and a Jordan $F$-algebra $J$, there is a dual pair
\[  \Aut(C) \times \Aut(J)  \longrightarrow \mathcal{E} \]
for an appropriate group $\mathcal{E}$. Indeed, this observation underlies what is known as the Tits' construction of exceptional Lie algebras.  For example, taking $C$ to be the octonion algebra $\mathbb{O}$, so that $\Aut(C)$ is the exceptional group $G_2$, one obtains a tower of dual pairs as $J$ varies, as tabulated below.
\vskip 5pt
 \begin{center}
 \begin{tabular}{|c|c|c|c|c|c|}
 \hline
 $\dim J$  &  3 & 6 & 9 & 15 & 27  \\
 \hline
 $\Aut(J)$ & $S_3$ & $\SO_3$ & $\PGL_3 \rtimes \Z/2\Z$ & $\PGSp_6$ & $F_4$  \\
 \hline
  $\mathcal{E}$ & $\PGSO_8 \rtimes S_3$ & $F_4$ & $E_6 \rtimes \Z/2\Z$ & $E_7$ & $E_8$ \\
  \hline
  \end{tabular}
 \end{center}
\vskip 5pt

\subsection{\bf Minimal representations}
What replaces the Weil representations? As noted before, the Weil representations are the smallest (in the sense of Gelfand-Kirillov dimension) infinite-dimensional representation of the metaplectic groups. It is their smallness that ensures that their restriction to dual pairs have manageable spectral decomposition. In the case of a general reductive group $\mathcal{E}$, one can thus seek the smallest infinite-dimensional representations of $\mathcal{E}$, the so-called minimal representations. 
The construction and classification of such  minimal representations have been carried out by Brylinski-Kostant  \cite{BryKostant94a,  BryKostant94b, BryKostant95} , Torasso \cite{Torasso97}, Kazhdan-Savin \cite{KazhdanSavin90}; see also \cite{GanSavin2005}. 
It turns out that for simply-connected groups $\mathcal{E}$'s of type $D$ or $E$, there is a unique minimal representation on the group $\mathcal{E}$, so that there is no need to consider nonlinear covers. In a sense, the situation is cleaner than that of the Weil representations.

\vskip 5pt

In the global setting, the global minimal representation (which is the restricted tensor product of the local ones) admits a (unique) embedding into the automorphic discrete spectrum $\mathcal{A}_2(\mathcal{E})$. This automorphic realization is not given by an ``averaging over rational points" type construction. Rather, it is constructed via residues of degenerate Eisenstein series \cite{GRS97a}. Ginzburg, Rallis and Soudry \cite{GRS97b} studied the global exceptional theta correspondence in many cases, and showed an analog of  the tower property for the tower of dual pairs  $\Aut(C) \times \Aut(J)$,  witth $C= \mathbb{O}$ fixed and as $J$ varies.  
\vskip 10pt

\subsection{\bf Howe duality}
With the main players in place, one can now consider exceptional theta correspondence. For example, using the minimal representation $\Pi$ of $\mathcal{E}$, we obtain a correspondence on ${\rm Irr}(\Aut(C)) \times {\rm Irr}(\Aut(J))$. For a representation $\pi \in {\rm Irr}(\Aut(C))$, one can define its big theta lift
\[  \Theta(\pi) = ( \Pi \otimes \pi^{\vee})_{\Aut(C)} \]
 which is a smooth representation of $\Aut(J)$. Unlike the situation of classical theta correspondence, the analog of the Howe duality theorem is not fully known. Part of this is due to the sporadic nature of the geometry of exceptional groups. In the archiemdean setting, there is the pioneering work of Huang-Pandzic-Savin \cite{HPS96} for the above tower of dual pairs with $Aut(J)$ compact, which implies a correspondence of infinitesimal characters \cite{Li99}. 
In the p-adic case, only recently have some cases been shown, namely for 
 \[ \Aut(\mathbb{O})  \times \Aut(J)  \quad \text{ for  $\dim J = 9$ or $15$,} \]
 i.e. for
 \[  G_2 \times \begin{cases}
\text{$ (\PGL_3 \rtimes \Z/2\Z)$ in $E_6 \rtimes \Z/2\Z$;} \\
\text{${\rm PD}^{\times}$ in $E_6^D$ (where $D$ is a cubic division algebra);} \\
\text{$\PU_3(K/F)$ in $E_6^K$;} \\
\text{$\PGSp_6$ in $E_7$.} \end{cases}\]
 These are contained in \cite{GanSavin2021} and \cite{BakicSavin2022}, and we shall give  a sketch of its proof below. Note however that there are known instances of the failure of the Howe duality property; see the forthcoming work of E. Karasiewicz for the dual pair $\PGL_2 \times F_4$ (taking $C$ to be the split quaternion algebra and $\dim J = 27$), which was first studied in \cite{Savin94}.
\vskip 5pt

\subsection{\bf Period transfers}   \label{SS:period transfers}
Recall that the key tools in the proof of the classical Howe duality theorem are the Jacquet modules of the Weil representation and  the doubling seesaw. In the exceptional setting,
the Jacquet modules of the minimal representation have been largely computed in the work of Savin and his collaborators \cite{MagaardSavin97, GrossSavin98, GanSavin2021}. There is however no analog of the doubling seesaw.
Hence, one needs to find a replacement.
\vskip 5pt

The point is to note that  the use of the doubling seesaw is a manifestation of a more general principle:
\[ \text{\it Theta correspondence frequently relates a period on one group with a period on the other}. \]
\vskip 5pt

Let us explain the precise meaning of this principle.
Suppose that $G \times H \rightarrow \mathcal{E}$ is a reductive dual pair, with a minimal representation  $\Omega$, over a local field $F$. 
Let $H'$ be  a subgroup of $H$ equipped with a character $\mu$. Then for $\pi \in {\rm Irr}(G)$, one has:
\begin{equation} \label{E:ss} 
 \Hom_{H'}(\Theta(\pi), \mu)  =   \Hom_{G \times H'} (\Omega,  \pi \boxtimes \mu) =  \Hom_G( \Omega_{H', \mu},  \pi). \end{equation}
Thus, if we could compute the twisted coinvariant  $\Omega_{H', \mu}$ of the minimal representation in some other way, we will obtain a description of $\Hom_{H'}(\Theta(\pi), \mu)$. In turns out that frequently, one has 
\begin{equation} \label{E:frequent}
  \Omega_{H', \mu}  \cong {\rm ind}_{G'}^G  \nu  \end{equation}
as $G$-modules (or contains the RHS as a submodule with ``spectrally negligible" quotient).  
Then an application of Frobenius reciprocity will give a natural isomorphism
\[  \Hom_{H'}(\Theta(\pi), \mu)    \cong  \Hom_{G'}( \pi^{\vee},  \nu^{-1}).  \]
Thus, we will obtain a statement of the form:
\[  \text{$\Theta(\pi)$ has nonzero $(H' ,\mu)$-period $\Longleftrightarrow$ $\pi^{\vee}$ has nonzero $\nu^{-1}$-period} \] 
at the local level.  This is the precise formula of the principle above. In the above argument, the key step is thus the independent computation of $\Omega_{H',\mu}$.
\vskip 5pt

 How does one establish the identity (\ref{E:frequent})? 
 When $H'$ is unipotent, this typically follows by a direct computation analogous to the computation of the Jacquet modules of $\Omega$. 
At the other extreme, when $H'$ is reductive,  the identity (\ref{E:frequent}) frequently arises via  a see-saw diagram: 
\[
 \xymatrix{
 G'  \ar@{-}[dr] \ar@{-}[d] & H   \ar@{-}[d] \\
  G \ar@{-}[ur] & H'}
\] 
where $H' \times G' \rightarrow \mathcal{E}$ is another dual pair, with $G \subset G'$ and $H' \subset H$. For $\pi \in {\rm Irr}(G)$, the identity (\ref{E:ss}) is precisely the seesaw identity:
\[  \Hom_{H'}(\Theta(\pi), \mu) \cong \Hom_{G \times H'}(\Omega, \pi \boxtimes \mu)  \cong \Hom_{G}(\Theta(\mu), \pi), \]
the additional information here being that $ \Omega_{H', \mu} = \Theta(\mu)$  has an action by a larger group $G'$.
 Typically, the $G'$-module  $\Theta(\mu)$ is a simpler representation,  such as a degenerate principal series representation, and the space $ \Hom_{G}(\Theta(\mu), \pi)$ of co-periods has a better chance of being understood. In the global setting, $\Theta(\mu)$ may be the special value of  an Eisenstein series, and its pairing with $\pi$ may be a global zeta integral which represents an automorphic L-function of $\pi$.

\vskip 5pt

Returning to the two key tools in the proof of Howe duality in the classical setting, we note that they are both instances of such period transfers:
\vskip 5pt

\begin{itemize}
\item The computation of the Jacquet modules of $\Theta(\pi)$, so that $H' =N$ is the nilpotent radical of a maximal parabolic subgroup $P = MN$ of $H$, and $\mu$ is the trivial character. 

\vskip 5pt

\item For the doubling seesaw, one is relating the $\O(V)^{\Delta}$-period on $\O(V) \times \O(V)$ with the $\Sp(W)^{\Delta}$-period on $\Sp(W) \times \Sp(W)$.
\end{itemize}
Many  other instances of such period transfers have been shown in the setting of classical theta correspondence, though we did not have the occasion to mention them in our earlier discussion.
 
\vskip 5pt

\subsection{\bf Period ping-pong} Let us give some examples of the period transfers relevant for the proof of Howe duality for, say, $G_2 \times PD^{\times}$ where $D$ is a cubic division algebra \cite{GanSavin2021}. In this case, since $PD^{\times}$ is compact, one can write:
\[  \Omega = \bigoplus_{\pi \in {\rm Irr}(PD^{\times})} \pi \otimes \Theta(\pi)  \]
with $\Theta(\pi)$ a semisimple unitarizable representation of $G_2$. Likewise $\Theta(\sigma)$ is semisimple for $\sigma \in {\rm Irr}(G_2)$, and is nonzero only when $\sigma$ is tempered.
\vskip 5pt
 
Let $P = M\cdot N$ be the so-called Heisenberg parabolic subgroup of $G_2$, so that $M \cong \GL_2$ and $N$ is a 5-dimensional Heisenberg group. The generic $M$-orbits of characters of $N$ are naturally parametrized by the isomorphism classes of \'etale cubic $F$-algebras $E$. On the other hand, the set of cubic field extensions of $F$   parametrize the maximal tori in $D^{\times}$, since every cubic field $E$ admits an embedding into $D$ uniquely up to conjugacy. Now for $\pi \in {\rm Irr}(PD^{\times})$ and $E$ an \'etale cubic $F$-algebra with associated $\psi_E \in \Hom(N, \C^{\times})/_M$, the period transfer identity gives:
\[  \Hom_N( \Theta(\pi), \psi_E)  \cong \Hom_{PD^{\times} \times N}(\Omega, \pi \otimes \psi_E) \cong \Hom_{PD^{\times}}( \Omega_{N, \psi_E}, \pi). \]
It turns out by a direct computation that one has
\[  \Omega_{N, \psi_E} \cong \begin{cases}
 {\rm ind}^{PD^{\times}}_{PE^{\times}} 1, \text{  if $E$ is a field;} \\
 0, \text{  if $E$ is not a field.}
 \end{cases} \]
 Hence one concludes (taking note that $\pi^{\vee} \cong \bar{\pi}$ is unitary) that:
\[  \Hom_N( \Theta(\pi), \psi_E)    \cong \begin{cases}
 \Hom_{PE^{\times}}(\pi, \C), \text{  if $E$ is a field;} \\
 0, \text{  if $E$ is not a field.}   \end{cases} \]
 In other words, the $(N, \psi_E)$-period of $\Theta(\pi)$ is related to the torus period of $\pi$ with respect to $T_E = PE^{\times}$.
 \vskip 5pt

 We now play a game of period ping pong by  turning the table around:  start with $\sigma \in {\rm Irr}(G_2)$ and  compute the torus period $\Hom_{PE^{\times}}(\Theta(\sigma), \C)$. This is a reductive period and there is an associated  seesaw diagram:
\[
 \xymatrix{
 \Spin_8^E  \ar@{-}[dr] \ar@{-}[d] & \Aut(D) = PD^{\times}
     \ar@{-}[d] \\
  G_2 \ar@{-}[ur] &   \Aut(i: E \rightarrow D) \cong PE^{\times}}
\]
where $\Spin_8^E$ is the quasi-split form of   $\Spin_8$ determined by $E$. 
The dual  pair $PE^{\times} \times \Spin_8^E$ here is not of the form $\Aut(C) \times \Aut(J)$ but is associated to a pair of so-called twisted composition algebras (a notion due to Springer).
In any case, the seesaw identity  gives
\[  \Hom_{PE^{\times}}(\Theta(\sigma), \C) \cong \Hom_{G_2}(\Theta(1), \sigma). \]
The reader who still recalls our sketch of the Howe duality theorem in the classical setting will not fail to notice the parallel here. As in the classical case, the representation $\Theta(1)$ is a quotient of a degenerate principal series $I_E$ of $\Spin_8^E$, so that
\[  \Hom_{PE^{\times}}(\Theta(\sigma), \C) \cong \Hom_{G_2}(\Theta(1), \sigma) \subset \Hom_{G_2}(I_E, \sigma). \]
In addition, as a $G_2$-module, a Mackey theory type argument gives:
\[   {\rm ind}_N^{G_2} \psi_E \subset  I_E \quad \text{  with a small  associated quotient.} \]
This means that one has:
\[  \Hom_{PE^{\times}}(\Theta(\sigma), \C)   \subset \Hom_{G_2}(I_E, \sigma) \approxeq \Hom_N( \sigma^{\vee}, \psi^{-1}_E). \]
In fact, the last approximate equality  can be shown to be an equality when $\sigma$ is tempered. 
Hence we have related the $PE^{\times}$-period of $\Theta(\sigma)$ with the $(N,\psi_E)$-period of $\sigma$ (at least for tempered $\sigma$). 
\vskip 5pt

\subsection{\bf Sketch of Howe duality} This period ping pong is the key ingredient in proving the Howe duality theorem in the exceptional case \cite{GanSavin2021}. Indeed, assembling the information above, one sees that if 
\[  \Hom_{PD^{\times} \times G_2}(\Omega, \pi \otimes \sigma) \ne 0, \]
then one has a chain of containments:
\[  \Hom_N(\sigma, \psi_E) \subset \Hom_N(\Theta(\pi),\psi_E) \cong \Hom_{PE^{\times}}(\pi, \C) \subset \Hom_N(\sigma, \psi_E) \]
which close up into a cycle. Since we know that $\Hom_{PE^{\times}}(\pi, \C) $ is finite-dimensional (because $\pi$ is), we conclude that this is a chain of equality. From this, and the fact that $\Hom_N(\sigma, \psi_E)$ is nonzero for some $E$,  it is not hard to argue that one has
\[  \Theta(\pi) = \sigma \quad \text{and} \quad \Theta(\sigma) = \pi. \]

\vskip 5pt
\subsection{\bf Applications}
We conclude our brief discussion of exceptional theta correspondence by mentioning two applications:
\vskip 5pt

\begin{itemize}
\item in the local nonarchimedean setting, the theta correspondence for $G_2 \times H$ highlighted above was crucially used in the recent proof of a local Langlands correspondence for $G_2$ \cite{GanSavin2022};

\vskip 5pt

\item in the global setting, exceptional theta correspondence has been used to construct many nontempered A-packets of $G_2$; see \cite{GanGurevich2005} for a discussion.
In particular, it gives the most interesting nontempered A-packet of $G_2$: the cubic unipotent A-packets, which exhibits unbounded cuspidal multiplicities in the automorphic discrete spectrum of $G_2$ \cite{GanGurevichJiang2002}.
\end{itemize}
\vskip 5pt
Let us elaborate on this global application, as it is very nice and simple. Over a number field $k$, we consider the dual pair (with $J = F^3$ and $C = \mathbb{O}$):
\[  \Aut( F^3) \times \Aut(\mathbb{O})  = S_3 \times G_2 \longrightarrow \mathcal{E} = \PGSO_8 \rtimes S_3, \]
where $S_3$ is the symmetric group on three letters. The local minimal representation $\Omega_v$ decomposes as
\[  \Omega_v = \bigoplus_{\eta_v \in {\rm Irr}(S_3)} \eta_v \otimes \pi_{\eta_v} \]
and it was shown by Huang-Magaard-Savin  \cite{HMS98} (in the nonarchimedean case) and Vogan \cite{Vogan94} (in the archimedean case) that each $\pi_{\eta_v}$ is irreducible unitary. 
The abstract global minimal representation thus decomposes as:
\[  \Omega = \bigotimes'_v \Omega_v \cong \bigoplus_{\eta \in {\rm Irr}(S_3(\A))} \eta \otimes \pi_{\eta}, \quad \text{(with $\pi_{\eta} = \otimes_v \pi_{\eta_v}$)}. \]
Now recall we have an automorphic realization 
\[  \iota:  \Omega \hookrightarrow \mathcal{A}_2(\PGSO_8) \]
given by residues of an appropriate Eisenstein series. By pulling back functions on $[\PGSO_8]$ to $[G_2]$, it was shown in \cite{GanGurevichJiang2002} that one obtains a $G_2$-equivariant map
\[  {\rm Rest} \circ \iota : \Omega \twoheadrightarrow \Omega^{S_3(k)} \hookrightarrow \mathcal{A}_2(G_2). \]
This gives a submodule of $\mathcal{A}_2(G_2)$ isomorphic to
\[  \Omega^{S_3(k)}   \cong \bigoplus_{\eta}  \eta^{S_3(k)}  \otimes \pi_{\eta}. \]
Note that by the Peter-Weyl theorem,
\[  \dim \eta^{S_3(k)}  = \dim \Hom_{S_3(k)} ( \eta, L^2(S_3(k) \backslash S_3(\A))) \]
is the multiplicity of $\eta$ in the space of automorphic forms of the finite group scheme $S_3$. So what we have here is the analog of Theorem \ref{T:Lis}. 
A simple character computation shows that the set of numbers $\{ \dim \eta^{S_3(k)} \}$ is unbounded.
 
\vskip 5pt

\section{\bf Theta Correspondence and the Relative Langlands Program}

So far, we have given a rather traditional treatment of the theory of theta correspondence. However, the discussion of period transfers in the previous section provides us the opportunity to connect with the relative Langlands program.
In this section, our goal is to explain how the  theta correspondence  fits into the framework of the relative Langlands program.
\vskip 5pt

The relative Langlands program is built upon the pioneering and fundamental work of Herve Jacquet, in his study of periods of automorphic forms and their relation with Langlands functoriality. In particular, he introduced the relative trace formula, which is the basic tool in the subject. Via comparison of relative trace formulas, Jacquet envisioned that one can relate automorphic periods on different groups, thereby understanding them in terms of the Langlands program. The relatively recent book \cite{SakeVenk2017} of Sakellaridis-Venkatesh supplies a unifying conceptual  framework for the study of such automorphic periods, casting the subject as a bona-fide extension of the classical Langlands program. 
\vskip 5pt

On the other hand, as we saw in the previous section, the theory of theta correspondence provides an effective tool for the transfer of periods between members of a dual pair. It is thus not surprising that the theta correspondence has a role to play in the relative Langlands program, if nothing else, as a complementary  tool  to the relative trace formula. As we shall see, the connection in fact runs deeper than its utility as a tool.  Rather, it can be subsumed as an instance of the relative Langlands program.

\vskip 5pt
\subsection{\bf Relative Langlands}
We briefly recall the formulation of the relative Langlands program, as covered in Beuzart-Plessis's lectures.
In its initial conception, as given in the book \cite{SakeVenk2017} of Sakellaridis-Venkatesh, the relative Langlands program is concerned with a spherical subgroup $H \subset G$, so that $X= H \backslash G$ is a spherical variety.
For simplicity, we shall assume that $X(F) = H(F) \backslash G(F)$ for the fields $F$ we are considering below. 
We also take a unitary character $\chi: H(F) \longrightarrow S^1$. 
Roughly, one is interested in the following problems (see also \cite{GanWan2019}):
\vskip 5pt

\begin{itemize}
\item ($L^2$-setting) For a local field $F$, determine the spectral decomposition of 
\[   L^2(X(F), \chi) = L^2((H(F),\chi) \backslash G(F)) ={\rm Ind}_{H(F)}^{G(F)} \chi \]
 as a unitary representation of $G(F)$. 
\vskip 5pt

\item (Smooth local setting) Still over a local field $F$, classify the set
\[  {\rm Irr}_X(G(F)) = \{ \pi \in {\rm Irr}(G(F)): \Hom_{H}(\pi, \chi) \ne 0 \}\]
of $(H, \chi)$-distinguished irreducible smooth representations of $G(F)$.
Moreover, determine $\dim \Hom_{H}(\pi, \chi)$ when it is nonzero, and produce a natural basis of this Hom space if possible.
Indeed, a special case of particular interest is when there is a multiplicity-at-most-one situation, where $\dim \Hom_H(\pi,\chi) \leq 1$.
\vskip 5pt
By Frobenius reciprocity, the above problem is equivalent to classifying irreducible submodules of the smooth $G(F)$-representation 
\[  C^{\infty}(H(F), \chi \backslash G(F)) = {\rm Ind}_{H(F)}^{G(F)} \chi, \]
and by duality and complex conjugation, this is in turn equivalent to classifying the irreducible quotients of the compactly-induced representation
\[  C^{\infty}_c(H(F) , \chi \backslash G(F)) = {\rm ind}_{H(F)}^{G(F)} \chi. \]
Rephrased in this way, one sees that this is just the smooth version of the $L^2$-problem formulated above. 
\vskip 5pt

\item (Global setting) Let  $\chi: H(\A) \longrightarrow S^1$ be an automorphic character of $H$ and consider the automorphic period integral
\[  \mathcal{P}_{H, \chi} : \mathcal{A}_{cusp}(G) \longrightarrow \C \]
defined by
\[  \mathcal{P}_{H,\chi}(f) = \int_{[H]} f(h) \cdot \overline{\chi(h)} \, dh. \]
Then we would like to classify the $(H,\chi)$-distinguished cuspidal representations, i.e. those  $\pi \subset \mathcal{A}_{cusp}(G)$ on which  $\mathcal{P}_{H,\chi}$ is nonzero.
Further, in the multiplicity-at-most-one setting, we would like to have a factorization of $|\mathcal{P}_{H,\chi}|^2$ as an Euler product of natural  local functionals.
\vskip 5pt

One may reformulate the above in a slightly different way, taking $\chi =1$ for simplicity. For $\phi \in  C^{\infty}_c(X(\A))$, one may form an $X$-theta series by
\[  \theta_X(\phi)(g) = \sum_{x \in X(F)}  (g \cdot \phi) (x), \]
so that one has a $G(\A)$-equivairant map
\[  \theta_X: C^{\infty}_c(X(\A)) \longrightarrow \mathcal{C}([G]). \]
The global problem of classifying $H$-distinguished $\pi$ is equivalent to classifying those $\pi$ which contribute to the spectral expansion of $\theta_X(\phi)$ (as $\phi$ varies), i.e. those $\pi$ such that 
\[   \langle \theta_X(\phi),  f \rangle_{[G]} \ne 0 \quad \text{ for some $f \in \pi$.} \]
 Indeed, putting in the definition of $\theta_X(\phi)$, one has
 \begin{align}
    \langle \theta_X(\phi),  f \rangle_{[G]} &=  \int_{[G]} \sum_{x \in X(F)}  \phi(g^{-1} \cdot x)  \cdot \overline{f(g)} \, dg \notag \\
    &= \int_{H(k) \backslash G(\A)} \phi(g^{-1} \cdot x) \cdot \overline{f(g)} \, dg \notag  \\
    &= \int_{X(\A)}  \phi (x)  \cdot \mathcal{P}_H(f)(x)  \, dx  \notag
    \end{align}
The nonvanishing of $\mathcal{P}_H(f)$ is thus equivalent to the existence of some $\phi$ such that the inner product of $\theta_X(\phi)$ and $f$ is nonzero. 
The problem of Euler factorization of the global $H$-period can also be expressed in terms of the $X$-theta series.
\end{itemize}

\vskip 10pt

\subsection{\bf  Proposed answers} \label{SS:proposed}
In \cite{SakeVenk2017}, Sakellaridis and Venkatesh proposed conjectural answers to the above problems, by introducing the dual group $X^{\vee}$ of $X$. More precisely,  they associated to $X= H \backslash G$
\vskip 5pt

\begin{itemize}
\item a (connected complex) dual group $X^{\vee}$;
\item a morphism
\[  \iota: X^{\vee} \times \SL_2(\C) \longrightarrow G^{\vee}, \]
well-determined up to $G^{\vee}$-conjugacy;
\item a $\frac{1}{2} \Z$-graded finite-dimensional representation $V_X$ of $X^{\vee}$.
\end{itemize}
If we let $G_X$ be a split reductive group with dual group $X^{\vee}$, then
the map $\iota$ gives rise to a map
\[ \iota_*: \Psi(G_X) := \{ \text{$A$-parameters of $G_X$} \} \longrightarrow \{ \text{$A$-parameters of $G$} \} =: \Psi(G), \] 
given by:
\[  \iota_*(\psi) = \iota \circ (\psi \times {\rm id}_{\SL_2}): L_F \times \SL_2(\C) \longrightarrow G^{\vee}. \] 
In particular, one obtains a lifting of tempered L-parameters of $G_X$ to $A$-parameters of $G$.
Using these data, they formulated conjectural answers to the problems above:
\vskip 5pt

\begin{itemize}
\item[(a)]  For the $L^2$-problem,  
one has the spectral decomposition
\[  L^2(X) \cong \int_{\Psi_{temp}(G_X)}   \Pi(\psi) \, d\mu(\psi) \]
where $d\mu(\psi)$ is the (natural) Plancherel measure on the space $\Psi_{temp}(G_X))$ of tempered A-parameters of $G_X(F)$ and $\Pi(\psi)$ is a unitary representation whose summands belong to the A-parameter $\iota_*(\psi)$. 
\vskip 5pt

 \vskip 5pt

\item[(b)]  in the smooth setting, the irreducible representations which occur as quotients of $C^{\infty}_c(X)$ are  (to a first approximation) those belonging to A-parameters which factor through $\iota$. 
\vskip 5pt

\item[(c)]  Globally,  the global period integral $|P_{H,\chi}|^2$, when restricted to $\Pi = \iota_*(\Sigma)$,  can be expressed as the product of a global constant and an Eulerian product of canonical local functionals.  The global constant is the L-value 
\[  \frac{L(1/2, \Sigma, V_X)}{L(1, \Sigma, {\rm Ad})},   \]
 whereas the local functionals are inherited from the spectral decomposition in (a).
 \vskip 5pt

The above statements may seem somewhat coarse, but they provide  a systematic nontrivial constraint on the spectral support of $L^2(X)$ or $C^{\infty}_c(X)$ in terms of a dual object $\iota: X^{\vee} \times \SL_2(\C)$, as the Langlands philosophy dictates.  The statement (c) also gives a systematic prediction on how and which L-values 
appear in automorphic periods.  Finally, the above framework  suggests a relative Langlands functoriality principle:
 \vskip 5pt
 
\item[(d)] Suppose that $X$  is a spherical $G$-variety and $Y$ is a spherical $H$-variety,  and there is a commutative diagram:
\[  \begin{CD}
  X^{\vee} \times \SL_2(\C) @>\iota_X>> G^{\vee}   \\
  @VVV   @VVV  \\
  Y^{\vee} \times \SL_2(\C)   @>\iota_Y>>  H^{\vee}.  
  \end{CD} \]
  Then there is a corresponding Langlands functorial lifitng from $X$-distinguished representations of $G$ to $Y$-distinguished representations of $H$. 
\end{itemize}

\vskip 5pt

\subsection{\bf An example} \label{SS:an example}
As an example (which was studied in \cite{GanWan2019}), let $V$ be a nondegenerate quadratic space of even dimensional $2n$ containing a nondegenerate line $L$, so that $V = L \oplus L^{\perp}$. Consider the hyperboloid 
\[  X = \O(V) / \O(L^{\perp}) \cong \O_{2n}/ \O_{2n-1} \]
which is the $\O(V)$-orbit of a nonzero vector  $v_L \in L$. This is a spherical variety (indeed a symmetric variety) whose associated dual data is
\[  \iota: X^{\vee}  \times \SL_2(\C)  = \SO_3(\C)  \times \SL_2(\C) \longrightarrow \SO_3(\C) \times \SO_{2n-3}(\C) \subset \O(V)^{\vee} = \SO_{2n}(\C) \]
where the map $\SL_2\C) \rightarrow \SO_{2n-3}(\C)$ is the irreducible representation of dimension $2n-3$. 
\vskip 5pt

Since $\SL_2^{\vee} = \SO_3(\C)$, the conjecture predicts that representaitons of $\O(V)$ distinguished by $\O(L^{\perp})$ (i.e. with nonzero $\O(L^{\perp})$-fixed vectors) are associated to A-parameters  built out of representations of $\SL_2$ and the map $\iota$.  Moreover, an example of an $\SL_2$- spherical variety $Y$ with $Y^{\vee} = \SO_3(\C)$ is the Whittaker variety $(N, \psi)\backslash \SL_2$, with associated representation the Gelfand-Graev module $L^2(N,\psi \backslash \SL_2)$.  Hence, one expects a relative functorial lifting from the spectrum of $(N,\psi) \backslash \SL_2$ to the spectrum of $X$, via the map $\iota$. 
\vskip 5pt

Incidentally, we have seen this map $\iota$ before, in the setting of Adams' Conjecture \ref{conj:adams}.  

\vskip 5pt

\subsection{\bf Where stands theta correspondonce?}
Now the theory of theta correspondence does not arise from a spherical variety and so does not really fit into the above framework. On the other hand, it is not too far from it either, since in both cases, one is interested in understanding the spectrum of a representation of $G(F)$: the Weil representation $\Omega$ in the case of theta correspondence and 
the representation $C^{\infty}_c(X(F))$ in the case of spherical varieties. So the natural question is:
\begin{itemize}
\item What sort of natural $G(F)$-modules should the relative Langlands program be concerned with?
\end{itemize}
A rough answer is that it should be concerned with natural $G(F)$-modules which are multipliciry-free, or at worst have finite multiplicities. This is certainly the initial motivation for singling out spherical varieties. From this point of view, it is clear how theta correspondence (in view of the Howe duality theorem) fits in. 
\vskip 5pt

 In a recent development, Ben-Zvi, Sakallaridis and Venkatesh proposed a broader framework for the relative Langlands program which encompasses both spherical varieties and theta correspondence. This is the framework of {\em quantization of (appropriate) Hamiltonian $G$-varieties}. 
  
 \vskip 5pt

\section{\bf Hamiltonian $G$-varieties and Quantizations}
As an interlude, we remind the reader some pertinent points about Hamiltonian varieties and quantization, a process that lies at the root of representation theory. Much of the discussion in this section will be heuristic in nature. For a more detailed discussion of these notions, we refer the reader to the article \cite{Vogan87} of Vogan.
\vskip 5pt

\subsection{\bf Symplectic manifolds}
In classical mechanics, the phase space of a classical system (i.e. the moduli of all possible states of the given system) is modelled by a symplectic manifold $(\mathcal{M}, \omega)$, where $\omega$ is a nondegenerate closed symplectic 2-form on $\mathcal{M}$. The symplectic form gives an identification
\[  \iota_{\omega} : T\mathcal{M} \cong  T^*\mathcal{M} \]
of the tangent and cotangent bundles of $\mathcal{M}$.  The space $\mathcal{C}(\mathcal{M}, \R)$ of (smooth) $\R$-valued  functions on  $M$ is called the space of observables of the system. Any $f \in \mathcal{C}(\mathcal{M}, \R)$ gives a 1-form $df$ and hence a vector field $\mathfrak{X}_f$ via $\iota_{\omega}$. 
The symplectic form $\omega$ induces  a Poisson bracket on $\mathcal{C}(\mathcal{M},\R)$ (namely  a Lie bracket which is a derivation in each variable) via
\[  \{ f_1, f_2 \} = \omega(\mathfrak{X}_{f_1}, \mathfrak{X}_{f_2}), \]
 making  $\mathcal{C}(\mathcal{M}, \R)$ a Poisson algebra. In general, a manifold $\mathcal{M}$ for which $\mathcal{C}(M)$ is equipped with a Poisson algebra structure is called a Poisson manifold. Geometrically, a Poisson manifold admits a foliation whose leaves are symplectic manifolds. Many constructs and results in symplectic geometry continue to hold in the setting of Poisson manifolds and our discussion will happen in this broader framework.
\vskip 5pt

\subsection{\bf Hamiltonian $G$-manifolds and moment maps}
 The symmetries of the classical system is represented by the (smooth) action of a symmetry group $G$ on the phase space $\mathcal{M}$ as
symplectomorphisms.   Differentiating the $G$-action gives a map
\[  X: \mathfrak{g} \longrightarrow \{ \text{Vector fields on $M$}\}. \]
 There is a notion of such a symplectic (or Poisson) $G$-manifold (or variety) being Hamiltonian. 
 Concretely, it means that it comes equipped with a $G$-equivariant moment map
\[  \mu : M \longrightarrow \mathfrak{g}^* \]
where $\mathfrak{g} = {\rm Lie}(G)$. This map is the mathematical formulation of the notion of momentum and it should satisfy the following property:
the dual map
\[  \mu^*: \mathfrak{g} = \mathfrak{g}^{**} \longrightarrow \mathcal{C}(\mathcal{M}, \R) \]
is a Lie algebra homomorphism such that 
\[   \mathfrak{X}_{\mu^*(z)} =   X_z \quad \text{   for $z \in \mathfrak{g}$.  } \]
The geometry of this moment map controls to some extent the geometry of the Hamiltonian $G$-variety $M$ and hence the properties of the classical mechanical system it represents. 
\vskip 5pt

\subsection{\bf Examples}
It will be pertinent to give some examples:
\vskip 5pt
\begin{itemize}
\item surely the simplest example of symplectic manifold is a symplectic vector space $(\mathbb{W}, \langle-,-\rangle)$. If $G \subset \Sp(\mathbb{W})$ is a subgroup, then $V$ is a (linear) Hamiltonian $G$-variety. The associated moment map
\[  \mu: \mathbb{W} \longrightarrow \mathfrak{g}^* \]
is given by
\[  \mu(w)(x) = \langle x \cdot  w, w \rangle, \quad \text{for $x \in \mathfrak{g}$.} \]
If $G= \Sp(\mathbb{W})$, then as a $G$-module, $\mathfrak{g}^* \cong {\rm Sym}^2(\mathbb{W})$; in this interpretation, the moment map is given as: $\mu(w) = w^2$.   
\vskip 5pt

\item The most basic example of a nonlinear symplectic manifold is the cotangent bundle $T^*(X)$ of a manifold $X$. If $X$ is a $G$-variety, then $G$  acts naturally on $T^*(X)$, making the latter a symplectic $G$-variety, for which the natural projection 
\[  \pi: T^*X \longrightarrow X \]
is $G$-equivariant.   In this case, the moment map is given by:
\[  \mu(x, \xi)( z) =  \langle d\pi_{(x,\xi)} \left(X_{z, (x,\xi)}\right), \xi \rangle, \quad \text{ for $x \in X$, $\xi \in T_x^*(X)$ and $z \in \mathfrak{g}$.} \]   

\vskip 5pt

\item In Lie theory, the most basic example of symplectic manifold is a coadjoint orbit $\mathcal{O} \subset \mathfrak{g}^*$.  For $\xi \in \mathcal{O}$ with stabilizer $G_{\xi} \subset G$, one can identify 
\[  T_{\xi}\mathcal{O}  = \mathfrak{g}/ \mathfrak{g}_{\xi}  \qquad \text{(where $\mathfrak{g}_{\xi} = {\rm Lie}(G_{\xi})$),} \] 
and the symplectic form on $T_{\xi}\mathcal{O}$ is given by:
\[  \omega_{\xi} (z_1, z_2) = \xi ( [ z_1, z_2]) \quad \text{   for $z_1, z_2 \in \mathfrak{g}$.}  \]
 The moment map is simply the natural inclusion $\mathcal{O} \hookrightarrow \mathfrak{g}^*$. 
Indeed, the partition of $\mathfrak{g}^*$ into $G$-orbits is a foliation with symplectic leaves, so that $\mathfrak{g}^*$ is a (linear) Poisson manifold, whose associated moment map is the identity map.  
\vskip 5pt

\noindent Having now realized that $\mathfrak{g}^*$ is a Poisson variety, we note that for a general Hamiltonian $G$-variety $\mathcal{M}$, the moment map $\mu: \mathcal{M} \rightarrow \mathfrak{g}^*$ is a Poisson algebra morphism.
\end{itemize}
These three examples will be the main ones we will consider. 
\vskip 5pt

\subsection{\bf Symplectic reduction}
Given a Hamiltonian $G$-variety $(\mathcal{M},\mu)$, because of the presence of symmetry through the $G$-action, it is natural to expect that one does not need the full symplectic manifold $(\mathcal{M}, \omega)$ to represent the classical system; perhaps one could cut down the number of parameters and use a smaller phase space. This process, which uses the moment map crucially,  is  a  basic construction in symplectic geometry known as  the Marsden-Weinstein symplectic (or Hamiltonian) reduction.  
\vskip 5pt

More precisely, let $\mathcal{O} \subset \mathfrak{g}^*$ be a $G$-orbit. Its preimage  
\[  \mu^{-1}(\mathcal{O}) = \mathcal{M} \times_{\mathfrak{g}^*} \mathcal{O}. \]
under the moment map $\mu$ inherits a $G$-action, and one can consider the (GIT) quotient 
\[   \mathcal{M}_{\mathcal{O}} :=  ( \mathcal{M} \times_{\mathfrak{g}^*} \mathcal{O}) //_G. \]
Under favourable conditions, $\mathcal{M}_{\mathcal{O}}$ inherits a nondegenerate symplectic form $\omega^{\flat}$  from $(\mathcal{M}, \omega)$, satisfying 
\[  p^*(\omega^{\flat}) =    \omega |_{\mu^{-1}(\mathcal{O})}, \]
with $p: \mathcal{M} \rightarrow \mathcal{M}_{\mathcal{O}}$ the natural projection,  so that $(\mathcal{M}_{\mathcal{O}}, \omega^{\flat})$ is a symplectic variety.  Indeed, one typically takes $\mathcal{O}$ to be the zero coadjoint orbit. 
\vskip 10pt

\subsection{\bf Quantization}
In passing from classical mechanics to quantum mechanics, one seeks to replace a classical system with a quantum system: this process is called quantization.
As one learns from a first course in quantum mechanics, this means replacing the state space $(\mathcal{M}, \omega)$ by a Hilbert space $V_{\mathcal{M}}$, or rather its projectivization $\mathbb{P}(V_{\mathcal{M}})$. The quantum state of a system is, for example, represented by a unit vector in $V_{\mathcal{M}}$. The classical observables in $\mathcal{C}(\mathcal{M},\R)$ are then replaced by Hermitian operators on $V$, whose eigenvalues give the quantum spectrum of the classical observable.  Further, if $G$ acts on $(\mathcal{M}, \omega)$, then its quantization $V_{\mathcal{M}}$ inherits a $G$-action as unitary operators. In other words, $V_{\mathcal{M}}$ becomes a unitary representation of $G$. 
\vskip 5pt

 As practised by the theoretical physicists of a century ago, this process of quantization is as much an art as a science, so we shall continue our discussion in the same spirit. The important takeaway for us is the following message:
 \[ \text{{\em Quantization replaces a symplectic $G$-variety  by a unitary representation of $G$.}}   \]
Moreover, the underlying Hilbert space of the unitary representation is typically obtained as follows. 
Recall that a Lagrangian submanifold   of $\mathcal{M}$ is a submanifold $\mathcal{Y}$ such that for each $y \in \mathcal{Y}$, $T_y(\mathcal{Y})$ is a maximal isotropic subspace of $T_y(\mathcal{M})$.
Let $\mathcal{X} \subset \mathcal{M}$ be a   submanifold with a Lagrangian fibration
\[   \mathcal{M} \longrightarrow  \mathcal{X}, \]
i.e. the fibers of this map  are Lagrangian submanifolds.   Then the Hilbert space in question should be the space of $L^2$-functions on $\mathcal{X}$, or more generally the space of $L^2$-sections of appropriate line bundles on $\mathcal{X}$.

\vskip 5pt

\subsection{\bf Quantization of examples}
With this message in mind, let us revisit the three examples of Hamiltonian $G$-variety above and examine what their quantizations should be:
\vskip 5pt

\begin{itemize}
\item For a symplectic vector space $\mathbb{W}$, consider a polarization $\mathbb{W} = \mathbb{X} \oplus \mathbb{Y}$. Then the projection $\mathbb{W} \longrightarrow \mathbb{Y}$ along $\mathbb{X}$ is a Lagrangian fibration. Hence
a quantization of $\mathbb{W}$ should be a unitary representation of $\Sp(\mathbb{W})$ realized on $L^2(\mathbb{Y})$.  We have encountered such a representation: it is a Weil representation. The fact that this quantization is a representation of $\Mp(\mathbb{W})$ and not $\Sp(\mathbb{W})$ is one of the deep mysteries of nature.  In any case, if $\Sp(W) \otimes \O(V)$ is a dual pair in $\Sp(\mathbb{W})$, we see that a quantization of the $\Sp(W) \times \O(V)$-Hamiltonian variety $\mathbb{W} = V \otimes W$ is a Weil representation $\Omega_{V,W,\psi}$. 
\vskip 5pt

\item Consider the cotangent bundle $T^*(X)$ of a $G$-variety $X$, which is a Hamiltonian $G$-variety. 
The natural projection $\pi: T^*X \longrightarrow X$ (with $X$ regarded as  the zero section)   is a Lagrangian fibration (which is stable $G$-equivariant in this case). Hence, the quantization of $T^*(X)$ should be the unitary representation $L^2(X)$ of $G$. If $X$ is a spherical variety, then this unitary representation is the one considered in the relative Langlands program. 
\vskip 5pt

\item Attaching an irreducible unitary representation $V_{\mathcal{O}}$ naturally to a coadjoint orbit  $\mathcal{O}$ is the dream of the orbit method. This dream was first realized when $G$ is a nilpotent group:  Kirillov  produced a bijection between the set of coadjoint orbits of $G$ and the unitary dual $\widehat{G}$. This was extended to Type 1 solvable groups by Auslander-Kostant.
For reductive Lie groups, this orbit philosophy does not work so cleanly, but people have not stopped trying to make it more precise;  for the latest developments concerning quantization of nilpotent coadjoint orbits of semisimple groups, the reader can consult the lectures of Mason-Brown. 
\vskip 5pt

To give an example, if $\mathcal{O} = \mathcal{O}_{min}$ is a minimal coadjoint orbit, then its quantization is a minimal representation of $G$.
\vskip 5pt

\item What about the quantization of the Poisson variety $\mathfrak{g}^*$ itself? Since $\mathfrak{g}^*$ is a family of symplectic $G$-varieties (the coadjoint orbits), the orbit method suggests that one should  think of $\mathfrak{g}^* /_G$ as the moduli space of irreducible unitary representations of $G$, i.e. as the unitary dual $\widehat{G}$. The fibration $\mathfrak{g}^* \rightarrow \mathfrak{g}^*/_G$ should then be quantized to the universal family of irreducible unitary representations. We do not know how to make more precise sense of this, but here is another thought. The GIT quotient $\mathfrak{g}^*//_G$ may be thought of as the generic unitary dual, in which case the fibration $\mathfrak{g}^* \rightarrow \mathfrak{g}^*//_G$ should quantize to the Gelfand-Graev module $L^2(N, \psi \backslash G)$ for a Whittaker datum $(N,\psi)$. This is probably getting too speculative to be useful. 
\vskip 5pt

On the other hand, one may consider the quantization of $\mathfrak{g}^*$ from the viewpoint of its ring of functions $C(\mathfrak{g}^*)$. 
Indeed, a common quantization procedure is that of deformation quantization, in which one seeks to deform the commutative Poisson algebra $C(\mathcal{M}, \C)$ to a noncommutative filtered $\C$-algebra. From this point of view, it was argued in the paper \cite{Rieffel} of Rieffel that the (deformation) quantization of $C(\mathfrak{g}^*)$ is the (reduced) group $C^*$-algebra $C^*(G)$. More precisely,  the Fourier transform gives an isomorphism of algebras
\[   (C(\mathfrak{g}^*), \bullet)  \cong  (C(\mathfrak{g}, \ast) \]     
where the multiplication on $C(\mathfrak{g}^*)$ is pointwise multiplication and that on $C(\mathfrak{g})$ is convolution. Via this isomorphism, one transports the Poisson algebra structure from $C(\mathfrak{g}^*)$ to $C(\mathfrak{g})$. Then the exponential map provides a way of formally deforming from the group level to the Lie algebra level.   
\end{itemize}
\vskip 10pt

 \subsection{\bf Quantization of moment map}
Given the important role played by the moment map in symplectic geometry, it is natural to ask: what is its quantization? 
In view of the difficulty we encountered above in deciding what the quantization of $\mathfrak{g}^*$ should be, this is not entirely clear.  But let us again take the viewpoint of its ring of functions. One has the pullback
\[  \mu^*: C(\mathfrak{g}^*) \longrightarrow C(\mathcal{M}, \C), \]
as well as  the pushforward
\[  \mu_{!} : C_c(\mathcal{M}, \C) \longrightarrow C_c(\mathfrak{g}^*). \]
Let us try to quantize these instead:
\vskip 5pt

\begin{itemize}
\item as mentioned before, $C(\mathfrak{g}^*)$ should quantize to the (reduced) group $C^*$-algebra $C^*(G)$;
\vskip 5pt

\item moreover, $C(\mathcal{M}, \R)$ should quantize to the Hermitian operators on $V_{\mathcal{M}}$, so that $C(\mathcal{M}, \C)$ should then give  the whole endomorphism algebra ${\rm End}(V_{\mathcal{M}})$ and the subspace $C_c(V_{\mathcal{M}}, \C)$ of compactly supported functions should correspond to the subalgebra 
${\rm End}_{{\rm fin}}(V_{\mathcal{M}})$ of finite rank operators.
\end{itemize}
Hence, one expects to give maps
\[  C^*(G)  \longrightarrow {\rm End}(V_{\mathcal{M}}) \]
and
\[  {\rm End}_{{\rm fin}}(V_{\mathcal{M}})  \longrightarrow C(G). \]
 The first (which quantizes $\mu^*$) is just the unitary representation $(\rho_{\mathcal{M}}, V_{\mathcal{M}})$ of $G$ giving rise to an action of $C^*(G)$ on $V_{\mathcal{M}}$. The second (which quantizes $\mu_{!}$)  is the formation of matrix coefficients. More precisely,. an element $T \in    {\rm End}_{{\rm fin}}(V_{\mathcal{M}}) $ gives rise to the function $f_T$ on $G$ defined by
 \[  f_T(g) = {\rm Tr}(\rho_{\mathcal{M}}(g) \circ T). \]
\vskip 5pt

\subsection{\bf Quantization of symplectic reduction}
We can now consider the quantization of symplectic reduction. Since the quantization of the moment map is obtained on the level of rings of functions, we shall consider symplectic reduction from this viewpoint:
\[  C(\mathcal{M}_{\mathcal{O}}) =  C(\mathcal{M} \times_{\mathfrak{g}^*} \mathcal{O})^G = \left(C(\mathcal{M}) \otimes_{C(\mathfrak{g}^*)} C(\mathcal{O}) \right)^G. \]
If $V_{\mathcal{M}_{\mathcal{O}}}$ is the Hilbert space quantizing $\mathcal{M}_{\mathcal{O}}$, then ${\rm End}(V_{\mathcal{M}_{\mathcal{O}}})$ must be the quantization of the RHS.
Now the quantization of the RHS should be given by the following heuristic computation (in which we have used the duality of vector spaces liberally):
\begin{align}
  \left({\rm End}(V_{\mathcal{M}}) \otimes_{C^*(G)} {\rm End}(V_{\mathcal{O}}) \right)^G  
 &= \left(   V_{\mathcal{M}}^* \otimes   (V_{\mathcal{M}} \otimes V_{\mathcal{O}}^* )_G \otimes V_{\mathcal{O}} \right)^G  \notag \\
&= {\rm End}_G \left(   V_{\mathcal{O}}^* \otimes V_{\mathcal{M}},   (V_{\mathcal{M}} \otimes V_{\mathcal{O}}^* )_G \right) \notag \\
&= {\rm End}((V_{\mathcal{M}} \otimes V_{\mathcal{O}}^* )_G). \notag 
\end{align}
Hence, we arrive at the conclusion: the quantization of the symplectic reduction $\mathcal{M}_{\mathcal{O}}$ is
\[  (V_{\mathcal{M}} \otimes V_{\mathcal{O}}^* )_G \]
This space is essentially the (dual of the) multiplicity space of the irreducible representation $V_{\mathcal{O}}$ in $V_{\mathcal{M}}$. 
This interpretation of  multiplicity spaces in terms of symplectic reduction was first observed in a series of papers of Guillemin-Sternberg \cite{GS82, GS84} and goes by the yoga ``quantization commutes with reduction".

\vskip 10pt

Here is a table summarizing our heuristic discussion above.
\vskip 5pt

\begin{center}
\begin{tabular}{|c|c|c|c|}
\hline 
& Classical &   Quantum &  \\
\hline  
Hamiltonian variety & $(\mathcal{M}, \omega)$&     $(\rho_{\mathcal{M}}, V_{\mathcal{M}})$ & unitary rep. \\
\hline 
Functions & $C(\mathcal{M}, \R) \subset C(\mathcal{M}, \C)$  & ${\rm Herm}(V_{\mathcal{M}})  \subset {\rm End}(V_{\mathcal{M}})$ & Operators \\
\hline 
Moment map & $\mu: M \rightarrow \mathfrak{g}^*$   &    & \\
\hline 
Pullback via $\mu$  & $\mu^*: C(\mathfrak{g}^*) \rightarrow C(\mathcal{M},\C)$  &  $\rho_{\mathcal{M}}: C^*(G) \rightarrow {\rm End}(V_{\mathcal{M}})$ & $C^*$-alg. module \\
\hline 
Pushforward via $\mu$ & $\mu_{!}: C_c(\mathcal{M}, \C) \rightarrow C_c(\mathfrak{g}^*)$  & ${\rm End}_{{\rm fin}}(V_{\mathcal{M}}) \rightarrow C(G)$ & matrix coeff. \\
\hline 
Symplectic reduction & $\left( \mathcal{M} \times_{\mathfrak{g}^*} \mathcal{O}\right) //_G$  & $(V_{\mathcal{M}} \otimes \pi_{\mathcal{O}})_G$    & multiplicity space  \\
\hline
\end{tabular}
\end{center}

\vskip 10pt

\subsection{\bf Theta correspondence as quantization of symplectic reduction}
  We have arrived at the (philosophical) conclusion that theta correspondence can be regarded as the quantization of symplectic reduction. Indeed, 
in the context of theta correspondence, one has $\mathcal{M} = V \otimes W$ with the Weil representation $\Omega$ as quantization. 
Picking a coadjoint orbit $\mathcal{O} \subset \mathfrak{o}(V)^*$ corresponds to picking  $\pi \in \widehat{\O(V)}$, and the symplectic reduction
\[  \mathcal{M}_{\mathcal{O}} = \mathcal{M}  \otimes_{\mathcal{O}} \mathfrak{o}(V)^* //_{\O(V)}   \]
amounts to extracting the multiplicity space of the $\pi$-isotypic component of the Weil representation $\Omega$. As one has the commuting action of $\Sp(W)$, $\mathcal{M}_{\mathcal{O}}$ is a symplectic $\Sp(W)$-variety, whose quantization is the multiplicity space   
\[ \Theta(\pi) = (\Omega \otimes \pi^{\vee})_{\O(V)}, \]
which is the big theta lift of $\pi$.
\vskip 5pt

In the setting of $C^*$-algebra, this point of view has been noted and pursued in the work of Landsman \cite{Landsman}. For a recent work on theta correspondence from the viewpoint of $C^*$-algebra,  see the paper of \cite{MeslandSengun}.

\vskip 10pt

\section{\bf An Extended Relative Langlands Program}  
After this long foray into Hamiltonian varieties and quantizations, let us now return to the relative Langlands program.
\vskip 5pt

\subsection{\bf Extended RLP}
Our discussion above shows that both the theory of theta correspondence and the spectral theory of spherical varieties can be subsumed in a single framework, built upon the quantization of appropriate Hamiltonian $G$-varieties. We stress that this point of view on the Weil representation is not new at all; as we recalled at the beginning of this article, it is how the oscillator representation emerged from quantum mechanics in the first place!
\vskip 5pt

In view of the above, we should now extend the scope of the relative Langlands program  to concern the spectral decomposition of  representations of $G$ (in the smooth or $L^2$-setting) arising from appropriate Hamiltonian $G$-varieties.  

\vskip 5pt

\subsection{\bf The dual Hamiltonian variety}
 Now a key tenet of the Langlands philosophy is that there should be a dual object which controls the spectral problem. In the case of a spherical variety $X$, or rather of its cotangent bundle $T^*(X)$ from our new point of view, we have seen that the dual object consists of the data 
 \[  (\iota: X^{\vee} \times \SL_2(\C) \rightarrow G^{\vee}, V_X). \] 
What should the dual data be in the setting of a more general Hamiltonian $G$-variety? 
\vskip 5pt

It turns out that the dual data above can be repackaged as a Hamiltonian $G^{\vee}$-variety. This is most easily described in the special case where the following simplifying conditions hold:
\vskip 5pt

\begin{itemize}
\item the map $\iota$ is trivial on $\SL_2(\C)$; 
\item the graded representation $V_X$ is concentrated in degree $0$. 
\end{itemize}
Then an observation due to Ben-Zvi, Sakellaridis and Venkatesh is that the $G^{\vee}$-variety
\[   G^{\vee} \times^{X^{\vee}} V_X  \]
has a natural symplectic structure, with respect to which it is a Hamiltonian $G^{\vee}$-variety (over $\C$). 
\vskip 5pt

The above observation thus leads one to the following  suggestive picture:
\begin{itemize}
\item There is a class of Hamiltonian varieties with reductive group action, for which an involutive  theory of duality exists. This theory of duality should associate to each Hamiltonian $G$-variety $(\mathcal{M}, \omega)$ over a local field $F$ a dual Hamiltonian $G^{\vee}$-variety $(\mathcal{M}^{\vee}, \omega^{\vee})$ over $F$.
\vskip 5pt

\item  If $V_{\mathcal{M}}$ is the quantization of $(\mathcal{M}, \omega)$ (so that $V_{\mathcal{M}}$ is a unitary representation of $G(F)$), then the spectrral decomposition of $V_{\mathcal{M}}$ is governed by the symplectic geometry of $\mathcal{M}^{\vee}$ (where we regard $\mathcal{M}^{\vee}$ as a complex variety). 
 \vskip 5pt
 
 \item The above duality theory should be symmetric in $\mathcal{M}$ and $\mathcal{M}^{\vee}$. Thus the spectral decomposition of a quantization of $\mathcal{M}^{\vee}$ is likewise governed by the symplectic geometry of $\mathcal{M}$.
\end{itemize}
\vskip 5pt
\noindent The reader will be justified in complaining about the vagueness in the above statements. 
Indeed, there is as yet no well-defined theory of duality as envisioned above and no precise expression of how the geometry of $\mathcal{M}^{\vee}$ should control the spectral properties of the quantization of $\mathcal{M}$ (though the paper of Ben-Zvi, Sakellaridis and Venkatesh should offer some clues). There are so far only a (large) number of examples for which one can test these expectations. 
\vskip 5pt\

A main takeaway from the above discussion, however, is that in the extended relative Langlands program, restriction problems come in pairs, corresponding to quantizations of a pair of Hamiltonion manifolds in duality with each other.   
\vskip 5pt

\subsection{\bf Revisiting theta correspondence}
We shall now reconsider theta correspondence with these expectations in mind. One particularly interesting question may be:  what is the restriction problem which is dual to theta correspondence?
\vskip 5pt

We can first try to understand results in theta correspondence from the viewpoint of the proposed answers in \S \ref{SS:proposed}. 
Consider for simplicity the case when 
\[  \dim W = 2n-2 \quad \text{and} \quad \dim V = 2n. \]
Hence, we are considering the quantization of the symplectic variety
\[  \mathcal{M} = V \otimes W. \]
Given \S \ref{SS:proposed} (especially item (a)) and the results of Sakellaridis \cite{Sake2017} on the $L^2$-theta correspondence given in Theorem \ref{T:L2 problem}, it is reasonable to surmise that 
the dual data associated to the theta correspondence for the dual pair $\Sp(W) \times \O(V)$ (with $\Sp(w)$ smaller)  consists of the natural diagonal inclusion:
\[  G_X^{\vee} = \Sp(W)^{\vee} = \SO_{2n-1}(\C)\longrightarrow \Sp(W)^{\vee} \times \O(V)^{\vee} = \SO_{2n-1}(\C) \times \SO_{2n}(\C). \]
Moreover, the factorization of global period in the setting of theta correspondence is simply the Rallis inner product formula. From Theorem \ref{T:rallis inner product}, one sees that the (graded) representation $V_X$ satisfies:
\[  L(1, \Sigma, V_X) = L(1,  \Sigma, {\rm std}) \cdot L(1, \Sigma, {\rm Ad}), \]
 suggesting  that as a representation of $\SO_{2n-1}(\C) = \Sp(W)^{\vee}$,
\[ V_X \cong  {\rm std} \oplus {\rm Ad} \cong \mathfrak{so}_{2n}. \] 
Having pinned down these dual data, the associated dual symplectic variety is
 \[    \mathcal{M}^{\vee} = (\SO_{2n}  \times \SO_{2n-1}) \times^{\SO_{2n-1}} \mathfrak{so}_{2n}  = T^*( \SO_{2n-1} \backslash (\SO_{2n} \times SO_{2n-1}). \]
The quantization of this cotangent bundle is the representation
\[  L^2( \SO_{2n-1} \backslash (\SO_{2n} \times SO_{2n-1})). \]
So what is the branching problem dual to theta correspondence?  It is the branching problem addressed by the Gross-Prasad conjecture \cite{GGP2012}!
\vskip 10pt

\section{\bf Theta Correspondence and Relative Functoriality}

We have now explained how theta correspondence finds its natural place in the extended relative Langlands program. Given the utility of theta correspondence in  period transfers 
mentioned in \S \ref{SS:period transfers}, it should not be surprising that theta correspondence can be used to demonstrate instances of relative functoriality and the conjecture of Sakellaridis-Venkatesh highlighted in \S \ref{SS:proposed}.  
\vskip 5pt

Typically, one would expect such relative functoriality to be shown by a comparison of relative trace formulas, which will entail a definition of a transfer map between spaces of relative orbital integrals (with an associated fundamental lemma for the spherical Hecke algebras). Such a relative trace formula comparison thus proceeds by a geometric comparison of relative orbital integrals and treats all representations together as a whole. The spectral transfer of individual (packets of) representations is only extracted at the very end, as a consequence of a relative character identity. In contrast, the theta correspondence is from the onset a spectral transfer of individual representations.  In this section, we offer some clues to how the two approaches can be reconciled, by explaining in a simple example how the machinery of theta correspondence can be used to recover some features of the relative trace formula approach.

\vskip 10pt

\subsection{\bf $\O(L^{\perp})$-period on $\O(V)$.}  \label{SS:period on O(V)}
We shall return to the example mentioned in \S \ref{SS:an example}, where we have 
\[  V = L \oplus L^{\perp} \quad \text{and} \quad X = \O(V) /\O(L^{\perp}) = \O(V) \cdot v_L \subset V, \]
with $\dim V =2n \geq 4$ and $\dim L =1$. Recall that we have the the dual data
\[  \iota: X^{\vee} \times \SL_2(\C) = \SO_3(\C) \times \SL_2(\C) \longrightarrow \SO_3(\C) \times \SO_{2n-3}(\C) \subset \O(V)^{\vee} =\SO_{2n}(\C). \] 
The discussion in \S \ref{SS:an example} shows that the spectrum of $X$ is described by A-parameters obtained as functorial lifts from the spectrum of the Whittaker variety $(N,\psi) \backslash \SL_2$ via the map $\iota$. We also noted that the map $\iota$ is precisely the one featuring in Adams' Conjecture \ref{conj:adams}. All these suggest that 
theta correspondence for $\SL_2 \times \O(V)$ can be used to prove this instance of relative functoriality.  To this end, it was shown in \cite{GanWan2019} that:
\vskip 5pt

\begin{itemize}
\item ($L^2$ problem) Given the Whittaker-Plancherel theorem 
\[ L^2(N,\psi\backslash \SL_2) \cong \int_{\widehat{\SL_2}} \dim \sigma_{N,\psi} \cdot  \sigma \, d\mu_{\SL_2}(\sigma)  \]
where $d\mu_{\SL_2}$ is the Harish-Chandra Plancherel measure of the unitary dual $\widehat{\SL_2}$, one has
\[  L^2(X) \cong \int_{\widehat{\SL_2}} \dim \sigma_{N,\psi} \cdot  \theta_{\psi}(\sigma) \, d\mu_{\SL_2}(\sigma)  \]
where $\theta_{\psi}(\sigma)$ is (the unitary completion of) the theta lift of (the underlying smooth representation of) $\sigma$. 
\vskip 5pt

\item (Smooth problem) The period transfer identity in this case is:
\[\Hom_N( \Theta(\pi), \psi) \cong   \Hom_{N \times \O(V)}(\Omega, \psi \otimes \pi) \cong  \Hom_{\O(V)}( \Omega_{N,\psi}, \pi),   \]
where $\Omega = C^{\infty}_c(V)$ is the Weil representation and $\pi \in {\rm Irr}(\O(V))$. A simple computation shows that the restriction map from $V$ to $X$ induces an isomorphism
\[  C^{\infty}_c(V)_{N,\psi}  \cong C^{\infty}_c(X) \cong {\rm ind}_{\O(L^{\perp})}^{\O(V)} 1, \] 
  so that
  \[  \Hom_{N}(\Theta(\pi),\psi) \cong \Hom_{\O(L^{\perp})}(\pi^{\vee}, \C)  = \Hom_{\O(L^{\perp})}(\pi, \C), \]
  and
  \[  {\rm Irr}_X(\O(V)) = \{ \pi \in {\rm Irr}(\O(V)): \text{$\Theta(\pi)$ is $\psi$-generic on $\SL_2$} \}.\] 
 One checks that $\Theta(\pi)$ is irreducible, except possibly when $\pi$ is the (small) theta lift of the trivial representation of $\SL_2$. 
 Hence,  one sees that theta correspondence gives an injection
 \[  {\rm Irr}_{N,\psi}(\SL_2) \hookrightarrow  {\rm Irr}_X(\O(V)) \]
 which is almost surjective: the only representation which is possibly missed  in the target is the theta lift of the trivial representation of $\SL_2$.  Moreover, for $\sigma \in {\rm Irr}_{N,\psi}(\SL_2)$, one has an isomorphism 
 \[  \Hom_N (\sigma, \psi) \cong \Hom_{\O(L^{\perp})}(\theta(\sigma), \C) \]
 of 1-dimensional spaces. For tempered $\sigma$, one can use the solution of the $L^2$-problem to normalize this isomorphism, so that it is well-defined up to $S^1$ (instead of $\C^{\times}$).  
 \end{itemize}
\vskip 5pt

\subsection{\bf Spaces of test functions}
Moreover,  using the Weil representation, \cite{GanWan2019} developed a theory of transfers of test functions in the style of the relative trace formula. The smooth Weil representation of $\SL_2 \times \O(V)$ is realized on the Schrodinger model $C^{\infty}_c(V)$.   Now we construct the following diagram:
    \begin{equation}  \label{E:diagram}
 \xymatrix{ &  C^{\infty}_c(V)
 \ar[dl]_p \ar[dr]^q & \\
 \mathcal{S}(N, \psi \backslash \SL_2)
  & &
 C^{\infty}_c(X)  
  }
\end{equation}
Here the map
  \[  p:  C^{\infty}_c(V) = \Omega^{\infty} \longrightarrow  C^{\infty}(N,\psi \backslash \SL_2)  \]
is  given by
  \[  p(\Phi) (g)  :=  (g \cdot \Phi)(v_L), \]
  where $v_L \in L$ is a fixed nonzero vector, and the map
   \[  q = {\rm rest}:  C^{\infty}_c(V) \longrightarrow C^{\infty}_c(X) \]
  is simply the restriction map, given by
  \[  q(\Phi)(h) = (h \cdot \Phi)(v_L) = \Phi(h^{-1} \cdot v_L) \quad \text{ for $h \in  \O(L^{\perp}) \backslash \O(V) = X$.} \]   

\vskip 5pt

We take note of some properties of these maps and their images:
\vskip 5pt

\begin{itemize}
\item   The map $p$ is $\O(L^{\perp})$-invariant and $\SL_2$-equivariant.  Let us set
 \begin{equation} \label{E:test}
    \mathcal{S}(N,\psi\backslash \SL_2)   := \text{image of $p$}, \end{equation}
 noting that it is a $\SL_2$-submodule. In fact, it is contained in the Harish-Chadra Schwarz space of $\SL_2$. 
 \vskip 5pt
 
\item   The map $q$ 
  \[  q = {\rm rest}:  C^{\infty}_c(V) \longrightarrow C^{\infty}_c(X) \]
  is surjective and, as noted before,  induces an $\O(V)$-equivariant isomorphism
     \[  q:  C^{\infty}_c(V)_{N,\psi} \cong C^{\infty}_c(X)\]
 \end{itemize}    
The spaces $ \mathcal{S}(N,\psi\backslash \SL_2)$ and $C^{\infty}_c(X)$ are our spaces of test functions on the two spherical varieties in question. 
\vskip 10pt

\subsection{\bf Transfers}
 We now make a definition: 
\vskip 5pt

\begin{defn}  \label{D:trans}
Say that $f \in \mathcal{S}(N, \psi \backslash \SL_2)$ and $ \phi  \in C^{\infty}_c(X) $ are in correspondence (or are transfers of each other) if there exists $\Phi \in C^{\infty}_c(V)$ such that $p(\Phi)  =f$ and $q(\Phi) = \phi$. 
\end{defn}
By the surjectivity of $p$ and $q$ onto the spaces of test functions, it is immediate that every $f \in \mathcal{S}(N, \psi \backslash \SL_2)$ has a transfer $\phi \in C^{\infty}_c(X)$ and vice versa.
\vskip 5pt

The transfer correspondence can in fact be refined to a map if one passes to the setting of ``relative orbital integrals".  More precisely, the composite map
\[ \begin{CD}
  C^{\infty}_c(V) @>q>> \mathcal{S}(N,\psi \backslash \SL_2) @>>> \mathcal{S}(N,\psi \backslash \SL_2)_{N,\psi} \end{CD} \]
  is  $(N,\psi)$-invariant (in addition to being $\O(L^{\perp})$-invariant) and  thus factors as
  \[  \begin{CD}
  C^{\infty}_c(V) @>p>>  C^{\infty}_c(V)_{N,\psi} \cong S(X) @>\tilde{t}>>  \mathcal{S}(N,\psi \backslash \SL_2)_{N,\psi}  \end{CD} \]
The map $\tilde{t}$, being $\O(L^{\perp})$-invariant, factors further to give:
\[  t :  C^{\infty}_c(X)_{\O(L^{\perp})} \longrightarrow  \mathcal{S}(N,\psi \backslash \SL_2)_{N,\psi} . \]
This map $t$ is the transfer map (on the level of coinvairant spaces).
\vskip 5pt
A relative orbital integral on $X$ or a relative character of a $X$-distinguished representation is a linear form on $C^{\infty}_c(X)_{\O(L^{\perp})}$. Likewise a (tempered) relative orbital integral on the Whittaker variety $(N,\psi)\backslash \SL_2$ or a (tempered) relative character of a $(N,\psi)$-generic representation of $\SL_2$ is a linear form on $\mathcal{S}(N,\psi \backslash \SL_2)_{N,\psi}$. Hence the transfer map $t$ allows one to pull back and compare relative orbital integrals or relative characters of these two spherical varieties. 

\vskip 5pt
 
 \subsection{\bf Basic function and fundamental lemmas}
A basic tenet of relative trace formulae comparison is that, having defined transfers, one needs to establish the relevant fundamental lemma (for the spherical Hecke algebras) in the unramified setting. In such an unramified setting,  as discussed in \S \ref{SS:unramified reps}, the Weil representation possesses a distinguished vector $\Phi_0 \in C^{\infty}_c(V)$, which is the characteristic function of a self-dual lattice in $V$. We may assume that $v_L$ is contained in this self-dual lattice.
 \begin{defn}  \label{D:basic}
    Set 
    \[  f_0 = p(\Phi_0) \quad \quad \text{and} \quad \quad \phi_0 = q(\Phi_0). \]
    We call these the basic functions in the relevant space of test functions.
    \end{defn}
Observe that while $\phi_0$ is compactly supported, $f_0$ is not. Moreover, it follows from definition that one has the following ``fundamental lemma":
    \vskip 5pt
     \begin{lemma}  \label{L:fund}
     The basic functions $f_0$ and $\phi_0$ correspond.
     \end{lemma}
In  Theorem \ref{T:unramified}(b), we saw that there is a surjective homomorphism of spherical Hecke algebras 
 \[  \lambda:   \mathcal{H}(\O(V), K_V)   \longrightarrow \mathcal{H}(\SL_2, K) \]
such that for $\phi \in  \mathcal{H}(\O(V), K_V) $, $\phi - \lambda(\phi)$ kills $\Phi_0$. From this, one deduces:
 \begin{lemma}  \label{L:fund2}
  For any $f \in \mathcal{H}(\O(V), K_V)$,  the element $f \cdot \phi_0 \in C^{\infty}_c(X)$  corresponds to the element $\lambda(f) \cdot f_0 \in \mathcal{S}(N, \psi \backslash \SL_2)$.
    \end{lemma}
\vskip 5pt

\subsection{\bf  Relative character identity}
The theory of transfers discussed above allows one to formulate the transfer of relative characters. Let us briefly recall this notion of relative character. Suppose that $X = H \backslash G$ and $\pi \in {\rm Irr}_X(G)$ is a unitary representation equipped with an invariant positive inner product (so an isomorphism $\pi^{\vee} \cong \overline{\pi}$).
Given
\[  \ell \in \Hom_{G}(\pi, C^{\infty}(X)), \]
one deduces by duality and complex conjugation the  $G$-equivariant map
\[  \overline{\ell^{\vee}} : C^{\infty}_c(X) \longrightarrow \overline{\pi^{\vee}} \cong \pi, \]
where the last isomorphism comes from the inner product. 
The associated relative character associated to the pair $(\pi, \ell)$ is the linear form 
\[ \mathcal{B}_{\pi, \ell} : C^{\infty}_c(X) \longrightarrow \C \]
defined by
\[  \mathcal{B}_{\pi, \ell}(\phi) =   \ell( \overline{\ell^{\vee}}(\phi) )(1). \]
This relative character factors as:
\[  \mathcal{B}_{\pi, \ell} : C^{\infty}_c(X) \longrightarrow C^{\infty}_c(X)_H  \longrightarrow \C. \]
So $\mathcal{B}_{\pi, \ell}$ depends on fixing a point $[1] \in X$ and $H$ is the stabilizer of that point. 
If $\ell$ is $X$-tempered (a notion we won't define here), then $\mathcal{B}_{\pi,\ell}$  extends to the larger Harish-Chandra Schwarz space of $X$.
\vskip 5pt

In particular, we have the notion of relative characters for the Whittaker variety $(N,\psi) \backslash \SL_2$ and $X = \O(L^{\perp}) \backslash \O(V)$.
In our discussion in \S \ref{SS:period on O(V)}, we have noted that theta correspondence furnishes an isomorphism
\[   f:  \Hom_N (\sigma, \psi) \cong \Hom_{\O(L^{\perp})}(\theta(\sigma), \C) \]
which is normalized up top $S^1$ if $\sigma$ is tempered.  In \cite{GanWan2019}, one finds the following relative character identity:

\vskip 5pt

\begin{thm}
For $\sigma \in {\rm Irr}(\SL_2)$ a tempered $(N,\psi)$-generic representation and $\ell \in \Hom_N (\sigma, \psi)$,
\[  \mathcal{B}_{\sigma, \ell} \circ t  =  \mathcal{B}_{\theta(\sigma), f(\ell)}. \]
\end{thm}
Note that though $f(\ell) \in  \Hom_{\O(L^{\perp})}(\theta(\sigma), \C) $ is only well-defined up to $S^1$, the relative character $ \mathcal{B}_{\theta(\sigma), f(\ell)} $ is well-defined, since the $S^1$-ambiguity cancels out due to the occurrence of  $f(\ell) \circ \overline{f(\ell)^{\vee}}$.

\subsection{\bf Other instances}
The framework of transfers of test fucntions discussed here is quite robust. Here are some other interesting instances:
\vskip 5pt

\begin{itemize}
\item we noted that the doubling seesaw effects a period transfer relating $\Sp(W)^{\Delta} \backslash (\Sp(W) \times \Sp(W))$ and $\O(V)^{\Delta} \backslash \O(V) \times \O(V)$. The relative characters in this setting are simply the usual Harish-Chandra characters. Hence, the above setup results in a character identity relating $\pi$ and $\theta(\pi)$, at least in the almost equal rank case. 
\vskip 5pt

\item in the setting of exceptional theta correspondence, the analog of the example we discussed is the dual pair 
\[  \PGL_2 \times \Aut(J)  = \PGL_2 \times F_4 \subset  E_7 \]
 where $J = H_3(\mathbb{O})$ is the 27-dimensional exceptional Jordan algebra. The minimal representation $\Omega$ is realized as a space of smooth functions on the cone $C$ of rank 1 elements in $J$. The set of elements in $C$ of trace $1$ is a homogeneous space $X$ for $\Aut(J)$ and the stabilizer of an element $x_L$ is isomorphic to the group $\Spin_9$,
 so that $X \cong \Spin_9 \backslash F_4$.  The exceptional theta correspondence for $\PGL_2 \times \Aut(J)$ should lead to the relative functoriality for
\[  (N, \psi) \backslash \PGL_2 \quad \text{and} \quad  X= \Spin_9 \backslash F_4. \]

\vskip 5pt

\item  The (similitude) exceptional theta correspondence for
\[  \PGL_2 ^3 \times \PGSO_8  \quad \text{  in $E_7$} \]
 should lead to the relative functoriality for
\[  \left( (N, \psi) \backslash \PGL_2 \right)^3  \quad \text{and} \quad  X = G_2 \backslash \PGSO_8. \]

\vskip 5pt

\item  The (similitude) theta correspondence for 
\[  \GL_3 \times \GL_3 \times \GL_3 \quad \text{in $E_6$} \]
should lead to the relative functoriality for
\[  (N, \psi) \backslash \GL_3 \quad \text{and} \quad \GL_3^{\Delta} \backslash (\GL_3 \times \GL_3). \]
 \end{itemize}
\vskip 10pt
We close by mentioning that in \cite{Sake2021}, Sakellaridis has developed a general theory of transfers for rank 1 spherical varieties.

\vskip 10pt
\section{\bf Automorphic Descent}

We have now spent a considerable number of pages on the theory of theta correspondence, mainly because it is the most developed among techniques of explicit construction of automoprhic forms.  It will be remiss, however, not to mention other methods of explicit construction. In this section, we will discuss the technique of {\em automorphic descent}, a theory initially developed by Ginzburg, Rallis and Soudry \cite{GRS99, GRS2011}, and subsequently extended by Dihua Jiang and Lei Zhang \cite{JiangZhang2020}. Automorphic descent
gives one a way of constructing  the backward functorial lifting from general linear groups to classical groups.  We begin by explaining the main underlying idea.
\vskip 5pt

\subsection{\bf The idea} 
Suppose that one has a subgroup
\[  H \ltimes N \subset G \]
where $N$ is unipotent and normalized by $H$. Hence, $H$ acts on the Pontryagin dual $\widehat{N} = \Hom(N, S^1)$. Given $\psi \in\widehat{N}$, let $H_{\psi}$ denote the stabilizer of $\psi$ in $H$. Now suppose that $\Pi$ is a representation of $G$. Then the coinvariant space $\Pi_{N, \psi}$ is naturally a representation of $H_{\psi}$. 
In this way, we have produced a representation of $H_{\psi}$ from a representation of the larger group $G$. Naturally, we are interested in the setting when the representation of $H_{\psi}$ so obtained is irreducible (or close to such). For this, it is intuitively clear that the input representation $\Pi$ should not be too big. 
\vskip 5pt

\subsection{\bf Automorphic descent}
Let us now consider the classical automorphic descent of Ginzburg-Rallis-Soudry.  We shall focus on the example of $\SO_{2n+1}$ and give a brief description of the descent construction. Suppose we are given a tempered L-parameter 
\[  \Psi =  \Pi_1 \oplus ...\oplus \Pi_r \] 
of $\SO_{2n+1}$, with $\Pi_i$ distinct cuspidal representations of some $\GL_{d_i}$'s of symplectic type, i.e. the exterior square L-function $L(s, \Pi_i, \wedge^2)$  has a pole at $s=1$ for each $i$. Moreover, $\sum_i d_i = 2n$. This global L-parameter gives rise to a near equivalence class of cuspidal representations in $\mathcal{A}_2(\SO_{2n+1})$.      
    \vskip 5pt
      
  We may regard $\Psi$ as an an automorphic representation
  \[  \Pi  = \Pi_1 \times ...\times \Pi_r  \]
  of $\GL_{2n}$ obtained by parabolic induction.  The purpose of automorphic descent is to construct from $\Pi$  the unique globally generic element in the near equivalence class of $\mathcal{A}_2(\SO_{2n+1})$ determined by $\Psi$. 
\vskip 5pt

  Regarding $\GL_{2n}$ as a Levi subgroup of a maximal parabolic subgroup $P$ of $\SO_{4n}$, one may consider 
  the parabolically induced representation 
    \[  I_P(s, \Pi) = {\rm Ind}_P^{\SO_{4n}} \Pi \cdot |\det|^s \quad \text{  of $\SO_{4n}$.} \]
  At the point $s=1/2$,  $I_P(1/2, \Pi)$ has a unique irreducible quotient  $J(1/2, \Pi)$. The Eisenstein $E(s, \Pi)$ associated to $I_P(s, \Pi)$ has a pole at $s = 1/2$ and  
  its (iterated) residue at $s=1/2$ gives an equivariant  map
  \[  I_P(1/2, \Pi) \twoheadrightarrow J(\Pi, 1/2) \hookrightarrow  \mathcal{A}_2(\SO_{4n}). \]
    \vskip 5pt
 
 Now in the context of the Gross-Prasad conjecture \cite{GGP2012}, one can consider the Bessel-Fourier coefficient of $J(\Pi, 1/2)$  with respect to $\SO_{2n+1}$. 
 More precisely, consider the parabolic subgroup of $\SO_{4n}$ given by 
 \[ M \ltimes N =  \left( \GL_1^{n-1} \times \SO_{2n+2} \right) \ltimes N    \]
 where $N$ is the unipotent radical. Let $\psi$ be a generic automorphic character of $N$, i.e. an element of $\widehat{N}$ whose $M$-orbit is maximal.  
 The space $\widehat{N}$ can be identified with the abelianization 
 \[  \overline{N}^{ab} (k) \cong  k^{n-2} \times V_{2n+2}, \]
 where $V_{2n+2}$ is the standard representation of $\SO_{2n+2}$.
 The action of the Levi factor $\GL_1^{n-1} \times \SO_{2n+2}$ is given by:
 \[  (t_1, t_2,..., t_{n-1};  h) :  (x_1, ...,x_{n-2} , v)  = \left( \frac{t_2}{t_1} x_1, \frac{t_3}{t_2}x_2,..., \frac{t_{n-1}}{t_{n-2}} x_{n-2}; t^{-1}_{n-1} h \cdot v \right). \]
 A generic character is represented by an element 
 \[  (x_1,...,x_{n-2}, v) \quad \text{with $x_i \ne 0$ for all $i$, and $v$ is non-isotropic.} \]
 Hence the stabilizer $M_{\psi}$ of such a generic $\psi$ is given by
 \[  M_{\psi} = \text{the stabilizer  in $\SO_{2n+2}$ of $v$} \cong \SO_{2n+1}. \]   
 In particular, one has the Bessel subgroup
 \[  M_{\psi} \cdot N = \SO_{2n+1} \cdot N \subset \SO_{4n}. \]
 
 \vskip 5pt
 
 Abstractly, one may now consider the twisted coinvariant space
 \[  J(\Pi, 1/2)_{N,\psi}, \quad \text{which is a representation of $\SO_{2n+1}$.} \]
 Globally, one considers the space of automorphic functions on $\SO_{2n+1}$ obtained by computing the $(N,\psi)$-Fourier coefficient of the automorphic forms in $J(\Pi, 1/2)$. 
 Namely, for $f \in J(\Pi, 1/2)$, one considers its Fourier coefficient
 \begin{equation}
   f_{N,\psi}(h) = \int_{[N]} \overline{\psi(n)} \cdot f(n h) \, dn,  \quad \text{$h \in M_{\psi}(\A) = \SO_{2n+1}(\A)$.}  \end{equation}
 The space
 \[  \mathcal{D}(\Pi):= \{ f_{N,\psi}:  f \in J(\Pi, 1/2) \}  \]
 is called a Bessel descent of $\Pi$ and gives  an automorphic representation of $\SO_{2n+1}$.  
 Ginzburg, Rallis and Soudry \cite{GRS2011} showed:
 \vskip 5pt
 
 \begin{thm}
 The Bessel descent $\mathcal{D}(\Pi)$  is nonzero irreducible and cuspidal. It is the irreducible globally generic cuspidal representation belonging to the near equivalence class determined by $\Psi$. In particular, the weak Langlands functorial lift of $\mathcal{D}(\Pi)$ to $\GL_{2n}$ is $\Pi$.  
 \end{thm}
 \vskip 5pt
 
 While we have described the automorphic descent in the context of $\SO_{2n+1}$, it should not surprise the reader that the analogous construction applies to 
 all classical groups (including metaplectic groups). For the symplectic and metaplectic groups, the automorphic descent will involve the consideration of Fourier-Jacobi models rather than the Bessel models. 
 \vskip 5pt
 
 Such explicit construction of the (globaly generic) cuspidal representations in generic A-packets  can be more useful in practice than a mere classification result. For example, in the work of Lapid and Mao \cite{LapidMao}, they exploited the automorphic descent to give a precise computation of the Whittaker-Fourier coefficients of globally generic cuspidal representations of metaplectic groups, in the style of the Ichino-Ikeda conjecture. 
 
 \vskip 5pt
 
 \subsection{\bf Motivation: a Rankin-Selberg integral}
 Let us explain how Ginzburg, Rallis and Sourdry are led to the above construction. For this, we assume for simplicity that $\Pi$ is cuspidal.
 In an attempt to identify the irreducible constituents of $\mathcal{D}(\Pi)$, it is natural (after showing its cuspidality) to consider the inner product of $\mathcal{D}(\Pi)$ with an arbitrary cuspidal representation $\sigma$ of  $\SO_{2n+1}$.  Now for $f \in J(1/2, \Pi)$, let $\tilde{f} \in I_P(1/2, \Pi)$ be a preimage of $f$. Then
 \begin{align}
   f_{N,\psi} (h) &= \int_{[N]} \overline{\psi(n)} \cdot  {\rm Res}_{s=1/2} E(s, \tilde{f})(n h) \, dn  \notag  \\
 &= {\rm Res}_{s=1/2} \left(  \int_{[N]} \overline{\psi(n)} \cdot  E(s, \tilde{f})(n h) \, dn  \right) \notag
 \end{align}
 Hence, one is led to consider the following inner product:
  \[    Z(s, \tilde{f}_s, \varphi) := \int_{[\SO_{2n+1}]}    \left(  \int_{[N]} \overline{\psi(n)} \cdot  E(s, \tilde{f})(n h) \, dn  \right) \cdot \overline{\varphi(h)} \, dh \]
where $\varphi \in \sigma \subset \mathcal{A}_{cusp}(\SO_{2n+1})$. It turns out that this integral is identically zero unless $\sigma$ is globally generic, in which case it  represents the L-function $L(s, \Pi\times \sigma)$ \cite{GRS98}, i.e.
 \[    Z(s, \tilde{f}_s, \varphi)  \simeq L(s +\frac{1}{2}, \Pi \times \sigma). \]
  In particular, 
  \[   {\rm Res}_{s=1/2}  Z(s, \tilde{f}_s, \varphi)   \simeq  {\rm Res}_{s = 1/2} L(s+ \frac{1}{2}, \Pi \times \sigma). \]
  Now observe:
  \vskip 5pt
  
  \begin{itemize}
  \item the LHS is the inner product of an element of $\mathcal{D}(\Pi)$ with an element of $\sigma$. 
  \vskip 5pt
  
  \item   the RHS  is nonzero if and only if $\Pi$ is isomorphic to the weak Langlands functorial lift of $\sigma$ to $\GL_{2n}$. 
  \end{itemize}
  Hence we see that the irreducible constituents of $\mathcal{D}(\pi)$ are precisely the globally generic cuspidal $\sigma$ which weakly lifts to $\Pi$.  This explains why the above theorem of Ginzburg, Rallis and Soudry holds.
  \vskip 5pt

 \vskip 5pt
 \subsection{\bf Twisted automorphic descent}  \label{SS:TAD}
 Extending the above work of Ginzburg-Rallis-Soudry, Jiang and Zhang \cite{JiangZhang2020} introduced 
 the twisted automorphic descent, which  offers the possibility of constructing the other members of the global L-packet  of $\SO_{2n+1}$ associated to $\Psi$. 
  In the above construction, Ginzburg-Rallis-Soudry showed that the automorphic descent of $J(\Pi, 1/2)$ is equal to the unique globally generic $\sigma_0$   in the L-packet.  Suppose that $\sigma$ is another  cuspidal representation of the global L-packet. How can one construct $\sigma$ by an analogous process as above? 
\vskip 5pt

The idea of Jiang-Zhang is as follows.  The cuspidal representation $\sigma$ will possess some nonzero Bessel coefficient with respect to a smaller even orthogonal group $\SO_{2m}$. Suppose that the Bessel coefficient of $\sigma$ down to $\SO_{2m}$ contains a cuspidal representation $\tau$ belonging to a generic  A-parameter of $\SO_{2m}$.   
The GGP conjecture implies that $\sigma$ is the unique representation in the global L-packet which supports this Bessel period with respect to $\tau$. 
This means that we have a way of distinguishing  the representation 
$\sigma$  from the other members of the L-packet   (just as one uses the Whittaker-Fourier coefficient to detect the representation $\sigma_0$). If this is how one is hoping to detect the representation $\sigma$ in any explicit descent-type construction, then it is only reasonable that one should somehow involve the auxiliary $\tau$ in the descent construction of $\sigma$. 
 \vskip 5pt
 
 Here then is the construction of Jiang-Zhang.  We first pick a $\tau$ as above.  The existence of such a $\tau$ is not obvious, especially with the requirement that  it should belong to a generic A-parameter, and   is in fact taken as a hypothesis in \cite{JiangZhang2020}. Hence the construction in \cite{JiangZhang2020} is so far conditional, though one would like to think that this hypothesis is not unreasonable. 
 In any case, under this hypothesis, one considers the standard module 
   \[  I(\Pi \otimes \tau, 1/2) =  \Pi \cdot |-|^{1/2}  \rtimes   \tau  \]
 of $\SO_{4n+2m}$ defined by induction from the parabolic subgroup with Levi factor $\GL_{2n} \times \SO_{2m}$. Let $J(\Pi \otimes \tau, 1/2)$ be the Langlands quotient of this standard module. It can be embedded into the automorphic discrete spectrum as an iterated residue of the corresponding Eisenstein series and one can then consider the  
Bessel descent  of $J(\Pi \otimes \tau, 1/2)$ down to $\SO_{2n+1}$.   Jiang-Zhang verified that one gets a cuspidal automorphic representation of $\SO_{2n+1}$, all of whose irreducible components $\sigma'$ belong to the relevant global L-packet (this involves a local unramified computation). They then computed the Bessel period  of $\sigma'  \otimes \tau$ and showed that it is nonzero, thus showing that $\sigma'  =  \sigma$.  
  \vskip 5pt
  
\subsection{\bf Another Rankin-Selberg integral}
As in the case of automorphic descent, the construction of Jiang and Zhang is informed by a Rankin-Selberg integral representing the L-function $L(s, \Pi \times  \sigma)$ for cuspidal representations $\sigma$ of $\SO_{2n+1}$. This Rankin-Selberg integral does not require  $\sigma$ to be globally generic but  to possess a Bessel-Fourier coefficient  with respect to a cuspidal representation $\tau$ with a generic A-parameter.  
\vskip 5pt

More precisely, for $\tilde{f}_s \in I(\Pi \otimes \tau, s)$ with associated Eisenstein series $E(s, \tilde{f})$ on $\SO_{4n+2m}$, one considers the zeta integral
 \[    Z(s, \tilde{f}_s, \varphi) := \int_{[\SO_{2n+1}]}    \left(  \int_{[N]} \overline{\psi(n)} \cdot  E(s, \tilde{f})(n h) \, dn  \right) \cdot \overline{\varphi(h)} \, dh \]
where $\varphi \in \sigma \subset \mathcal{A}_{cusp}(\SO_{2n+1})$. Here, $N$ is the unipotent radical of parabolic subgroup of $\SO_{4n+2m}$ with Levi factor
\[  M= \GL_1^{n+m-1}  \times \SO_{2n+2}, \] 
and $\psi$ is a generic character of $N$, whose stabilizer in $M$ is $\SO_{2n+1}$.  This zeta integral vanishes identically unless $\sigma$ has nonzero Bessel-Fourier coefficient with respect to $\tau$, in which case it represents the standard L-function of $\Pi \otimes \sigma$ \cite{JiangZhang2014}:
\[    Z(s, \tilde.{f}_s, \varphi) \simeq  L(s+ \frac{1}{2} , \Pi \times \sigma). \]
 On taking residue at $s = 1/2$, this gives
 \[  \langle \mathcal{D}_{\tau}(\Pi) , \sigma \rangle \ne 0 \]
 if and only if 
 \[   \text{$\sigma$ has nonzero  $\tau$-Bessel coefficient and ${\rm Res}_{s = 1/2} \, L(s + \frac{1}{2}, \Pi \times \sigma) \ne 0$}.  \]
This explains why  $\mathcal{D}_{\tau}(\Pi)$ captures precisely the desired representation $\sigma$ in \S \ref{SS:TAD}. 
\vskip 5pt

\subsection{\bf Relation with nontempered GGP} 
As one can see from the above discussion, the automorphic descent is based on the consideration of Bessel or Fourier-Jacobi coefficients of certain residual Eisenstein series. 
Thus it is not surprising that it is related to the GGP conjecture. While the classical GGP conjecture concerns the Bessel or Fourier-Jacobi coefficients of 
automorphic representations belonging to generic A-parameters, the residual Eisenstein series relevant here belongs to nongeneric A-parameters. There is however now an nontempered extension of the GGP conjecture which deals with nongeneric A-parameters. The automorphic descent  described above can be explained neatly in terms of this nontempered GGP conjecture; see \cite{GGP2020}.

\vskip 10pt
\section{\bf Generalized Doubling and Double Descent}
In this final section, we discuss a recent breakthrough which goes beyond the doubling zeta integral and automorphic descent. 
\vskip 5pt

Recall that the doubling zeta integral of Piatetski-Shapiro gives an integral representation of the standard L-function of any cuspidal representation of  $G \times \GL_1$, where $G$ is a classical group. Moreover, this doubling zeta integral plays an important role in the theory of theta correspondence. 
In a recent work of Cai-Friedberg-Ginzburg-Kaplan \cite{CFGK2019}, a generalized doubling zeta integral was introduced which represents the standard L-function for $G \times \GL_k$ for any $k \geq 1$. This beautiful construction has some profound implications which we will discuss below.

\vskip 5pt

\subsection{\bf Generalized doubling}
Let us explain the generalized doubling zeta integral in the setting of the symplectic group $\Sp(W) = \Sp_{2n}$.  Recall first the setting of the classical doubling zeta integral.
One considers the doubled space
\[  \mathbb{W} = W + W^- \]
which contains the Lagrangian subspace $W^{\Delta}$, with associated Siegel parabolic subgroup $P(W^{\Delta})$. One considers the corresponding degenerate principal series representation 
\[  I_{P(W^{\Delta})}(s, \chi) = {\rm Ind}_{P(W^{\Delta})}^{\Sp(\mathbb{W})}  \chi \cdot |\det|^s \]
and the associated Siegel Eisenstein series 
\[  E(s, \Phi_s) \quad \text{  for $\Phi_s \in  I_{P(W^{\Delta})}(s, \chi)$.} \]
Then the global doubling zeta integral takes the form
\[  Z(s, \Phi_s, f_1, f_2) := \int_{[\Sp(W) \times \Sp(W^-)] } E(s,\Phi_s,g) \cdot \overline{f_1(g_1)} \cdot f_2(g_2)  \, dg_1 \, dg_2 \]
for $f_1$, $f_2$ belonging to a cuspidal representation  $\sigma$ of $\Sp(W)$.  As mentioned, this integral represents the standard L-function of $\sigma \otimes \chi$:
\[  Z(s, \Phi_s, f_1, f_2) \simeq L(s +\frac{1}{2}, \sigma \times \chi).\]
 \vskip 5pt
 
 For the generalized doubling, one is seeking to represent the L-function $L(s, \sigma \times \tau)$ for a cuspidal representation $\tau$ of $\GL_k$. Not surprisingly, one needs to involve $\tau$ (which plays the role of $\chi$ above) in the construction. Instead of the Siegel principal series  induced from $\chi$ and the associated Siegel Eisenstein series, one considers principal series representation induced from a Speh representation built out of $\tau$. 
 \vskip 5pt
 More precisely:
 \vskip 5pt

 \begin{itemize}
\item[(a)]  Given a cuspidal representation $\tau$ of $\GL_k$, consider the parabolically induced representation
\[  \tau |-|^{\frac{2n-1}{2}} \times \tau |-|^{\frac{2n-3}{2}} \times......\times \tau |-|^{-\frac{2n-1}{2}}  \quad \text{  of $\GL_{2nk}$.}\]
 This has a unique irreducible quotient $\Delta(\tau, 2n)$ which is known as a Speh representation.
 By Moeglin-Waldspurger, one knows that $\Delta(\tau, 2n)$ admits an embedding into the automorphic discrete spectrum of $GL_{2nk}$ as the iterated residue of the corresponding Eisenstein series.
 \vskip 5pt
 
 \item[(b)]  Consider now the parabolically induced representation
\[  I ( s, \Delta(\tau,2n)) = {\rm Ind}_{P_{2nk}}^{\Sp_{4nk}} \Delta(\tau,2n) \cdot |\det|^s \]
where $P_{2nk}$ is the Siegel parabolic subgroup of $\Sp_{4nk}$. For $\Phi_s \in I ( s, \Delta(\tau,2n))$, let $\mathcal{E}(s, \Phi_s)$ be the corresponding Eisenstein series.
\vskip 5pt

\item[(c)]  Now we are going to define a specific Fourier coefficient of the Eisenstein series $\mathcal{E}(s, \Phi_s)$.  Let $Q  = L \cdot U \subset \Sp_{4nk}$ be the parabolic subgroup with Levi factor
\[  L = \GL_{2n} ^{k-1} \times \Sp_{4n} \]
and unipotent radical $U$.  More conceptually, $Q$ is the stabilizer of a flag
\[  X_1 \subset (X_1\oplus X_2) \subset...\subset (X_1 \oplus....\oplus X_{k-1})   \]
of isotropic spaces with $\dim X_i = 2n$, and 
\[  L = \GL(X_1) \times....\times \GL(X_{k-1}) \times \Sp(W_{4n}),  \]
where
\[  W_{4n} = (X_1 \oplus....\oplus X_{k-1})^{\perp} / (X_1 \oplus....\oplus X_{k-1}) \]
is a nondegenerate symplectic space of dimension $4n$. 
Then
\[  U^{ab} \cong   
\Hom(X_2, X_1) \times ...\Hom(X_{k-1}, X_{k-2}) \times \Hom(W_{4n}, X_{k-1} ) \]
with the natural action of $L$, and 
\[  \Hom(U, \G_a) \cong \Hom(X_2, X_1) \times ...\times \Hom(X_{k-1}, X_{k-2}) \times \Hom(X_{k-1}, W_{4n}). \]
Now a generic element of $\Hom(U, \G_a)$ has the form  
\[  \mathfrak{X} = ( A_1,...., A_{k-2} ;  T)  \]
such that
\begin{itemize}
\item  $A_i \in {\rm Isom}(X_{i+1}, X_i)$;
\item Image of $T$ is a $2n$-dimensional nondegenerate subspace of $W_{4n}$. 
\end{itemize} 
\vskip 5pt
For $\mathfrak{X}$ as above,  its stabilizer in $L$ is given by
\[    L_{\mathfrak{X}}  \cong \Sp({\rm Im}(T)) \times \Sp({\rm Im}(T)^{\perp}) \cong \Sp_{2n} \times \Sp_{2n} \]  

\vskip 5pt

\item[(d)] The element $\mathfrak{X}$ gives rise to a generic unitary character $\psi_{\mathfrak{X}}$ of $\U(\A)$ trivial on $\U(k)$ (by composition with a nontrivial additive character $\psi$ of $k \backslash \A$). Hence, we may consider the Fourier coefficient of an automorphic form on $\Sp_{4nk}$ with respect to $(U, \psi_{\mathfrak{X}})$.
 In particular, consider the $(U, \psi_{\mathfrak{X}})$-Fourier coefficient of the Eisenstein series $\mathcal{E}(s, \Phi_s)$ from (b):
\[  \mathcal{E}(s, \Phi_s)_{U, \psi_{\mathfrak{X}}} (g) = \int_{[U]} \mathcal{E}(s,\Phi_s, ng) \cdot \overline{\psi_{\mathfrak{X}}(n)} \, dn. \]
This defines an automorphic form on the stabilizer $L_{\mathfrak{X}} \cong \Sp_{2n} \times \Sp_{2n}$. 

\vskip 5pt
\item[(e)] 
For a cuspidal representation $\sigma$ of $\Sp_{2n}$, we may thus form the generalized doubling zeta integral
\[  \mathcal{Z}(s, \Phi, f_1, f_2) :=  \int_{[\Sp_{2n} \times \Sp_{2n}]}  \mathcal{E}(s, \Phi_s)_{U, \psi_{\mathfrak{X}}} (g_1, g_2)   \cdot \overline{f_1(g_1)} \cdot f_2(g_2) \, dg_1 \, dg_2,\]
for $f_1$, $f_2 \in \sigma$. 
\end{itemize}
The main result of \cite{CFGK2019} (see also \cite{Cai2021})  is that 
\[   \mathcal{Z}(s, \Phi, f_1, f_2) \simeq  L(s+\frac{1}{2}, \sigma \times \tau). \]
A crucial point here is that there is no condition (such as genericity condition) on the cuspidal representation $\sigma$ of $\Sp_{2n}$. 

\vskip 5pt

An interesting aspect of this generalized doubling zeta integral is that it can be analogously carried out for nonlinear covers of $\Sp_{2n}$ (and more generally of classical groups). 
For this, see \cite{Kaplan2019} and \cite{Cai2022}.
\vskip 10pt

\subsection{\bf Functorial lifting}
We will finish by giving two applications of the generalized doubling zeta integral. The first is that it leads to the development of a local theory of $\gamma$-factors, $L$-factors and $\epsilon$-factors for representations of $G \times \GL_k$ for classical groups $G$, in the style of what \cite{LapidRallis2005}  did with the classical doubling zeta integral. This has been carried out by Cai-Friedberg-Kaplan \cite{CFGK2022a}. 
\vskip 5pt

This theory of local factors allows one to have a good handle of the analytic properties of the L-functions $L(s, \sigma \times \tau)$. As a spectacular application, Cai, Friedberg and Kaplan have shown the following theorem in \cite{CFGK2022b}:
\vskip 5pt

\begin{thm} 
Every cuspidal representation of a classical group $G$ has a weak functorial lift  to the relevant general linear group.
\end{thm}

\vskip 5pt
This is achieved by an application of  the Converse Theorem of Cogdell and Piatetski-Shapiro, in the same way as \cite{CKPSS} did for globally generic cuspidal representations of classical groups.   It will be interesting to see if one can obtain a refined classification of the fibers of this functorial lifting and a determination of its image, in the style of Arthur's conjecture.  This would give an alternative approach to the results of Arthur \cite{Arthur2012}.

\vskip 5pt

\subsection{\bf Double descent}
The second application, due to Ginzburg and Soudry \cite{GinzburgSoudry2021}, is a  descent construction which works unconditionally (i.e. without genericity hypothesis as in the original automorphic descent or the assumption on Bessel coefficients in twisted automorphic descent).  They call their construction the {\em double descent}, since it is suggested by the generalized doubling zeta integral. 
\vskip 5pt

Let us illustrate with the example of $G = \Sp_{2n}$.  In this case, the functorial lifting is given by:
\[  \mathcal{A}_{cusp}(\Sp_{2n}) \longrightarrow \mathcal{A}(\GL_{2n+1}). \]
Suppose now that one has a cuspidal representation $\Pi$ of $\GL_{2n+1}$ whose symmetric square L-function $L(s, \Pi, {\rm Sym^2})$ has a pole at $s = 1$. 
Then such a $\Pi$ is known to be in the image of the above Langlands functorial lifting. The double descent constructs the whole preimage of $\Pi$ under this  functorial lifting. 
\vskip 5pt
 
 Let us consider the generalized doubling zeta integral, taking 
 \[  k = 2n+1 \quad \text{and} \quad \tau = \Pi. \]
Then one has:
\[   \int_{[\Sp_{2n} \times \Sp_{2n}]}  \mathcal{E}(s, \Phi_s)_{U, \psi_{\mathfrak{X}}} (g_1, g_2)   \cdot \overline{f_1(g_1)} \cdot f_2(g_2)\, \, dg_1 \, dg_2 \simeq L(s +\frac{1}{2}, \sigma \times \Pi). \] 
Now by the above theorem, $\sigma$ has a weak lift $\Pi_{\sigma}$ on $\GL_{2n+1}$ and $\Pi_{\sigma}$ is self-dual. Hence,  the RHS of the above identity is
\[  L(s+\frac{1}{2}, \sigma \times \Pi) = L(s+\frac{1}{2} , \Pi_{\sigma}^{\vee} \times \Pi), \]
which   has a pole at $s=1/2$ if and only if $\Pi = \Pi_{\sigma}$. Taking residue at $s = 1/2$, we thus conclude that 
\[  \int_{[\Sp_{2n} \times \Sp_{2n}]}  {\rm Res}_{s= 1/2} \mathcal{E}(s,\Phi_s)_{U, \psi_{\mathfrak{X}}} (g_1, g_2) \cdot \overline{f_1(g_1)} \cdot f_2(g_2) \, dg_1 \, dg_2 \ne 0 \]
if and only if $\Pi$ is the weak lift of $\sigma$. 
\vskip 5pt

Motivated by the above discussion, one now defines the double descent of $\Pi$ by
\[  \mathcal{DD}(\Pi) :=  \text{the span of ${\rm Res}_{s= 1/2}  \mathcal{E}(s,\Phi_s)_{U, \psi_{\mathfrak{X}}}$} \subset \mathcal{A}(\Sp_{2n} \times \Sp_{2n}) \]
Ginzburg and Soudry showed \cite{GinzburgSoudry2021}:
\vskip 5pt

\begin{thm}
For a cuspidal representation $\Pi$ of $\GL_{2n+1}$ such that $L(s, \Pi, {\rm Sym}^2)$ has a pole at $s = 1$, its  double descent $\mathcal{DD}(\pi)$ is contained in the space of cusp forms of $\Sp_{2n} \times \Sp_{2n}$. Moreover, for any irreducible cuspidal representation $\sigma_1 \otimes \sigma_2^{\vee}$ of $\Sp_{2n} \times \Sp_{2n}$, one has
\[  \langle \mathcal{DD}(\Pi),  \sigma_1 \otimes \sigma_2^{\vee} \rangle \ne 0  \]
if and only if
\[ \sigma_1 \cong \sigma_2 \quad \text{and $\Pi$ is the weak lift of $\sigma_1 \cong \sigma_2$.} \]
\end{thm}
Thus, this theorem of Ginzburg-Soudry provides a construction of the full near equivalence class of $\mathcal{A}_2(\Sp_{2n})$ associated to the generic A-parameter determined by $\Pi$. A/gain, it will be interesting to see if one can use this as a starting point to obtain a fine classification of this near equivalence class in the style of Arthur's conjecture. 
On the other hand, taking the fine classification provided Arthur's results in  \cite{Arthur2012} as given, it will be interesting to see if this explicit construction of the near equivalence class can be used to establish the Ichino-Ikeda type conjecture in the context of the Gross-Prasad conjecture.

\end{document}